\documentclass[11pt]{article}
\usepackage{graphicx}
\usepackage{graphics}
\usepackage{amssymb,amsmath,amsthm,epsfig}
\newcommand{\J}{P^{(\alpha,\beta)}}
\newcommand{\wJ}{\widetilde P^{(\alpha,\beta)}}
\newtheorem{theorem}{Theorem}
\newtheorem{lemma}{Lemma}
\newtheorem{corollary}{Corollary}

\newtheorem{proposition}{Proposition}

\newtheorem{remark}{Remark}
\newcommand{\R}{\mathbb R}

\newcounter{reh}
\newcounter{rek}
\begin{document}

\begin{center}
{\large {\bf Random matrices based schemes for stable and robust nonparametric
and functional regression estimators.}}\\
\vskip 1cm Asma Ben Saber$^a$ and Abderrazek Karoui$^a$ \footnote[1]{Corresponding author: Abderrazek Karoui\\
Email addresses:  Abderrazek.Karoui@fsb.rnu.tn (A. Karoui), asmabensaber@gmail.com  (A. Ben Saber)\\
This work was supported in part by the  
 DGRST  research grant  LR21ES10 and the PHC-Utique research project 20G1503.}  
\end{center}
\vskip 0.05cm {\small

\noindent $^a$ University of Carthage,
Department of Mathematics, Faculty of Sciences of Bizerte, Tunisia.
}\\

\noindent{\bf Abstract}--- In the first part of  this work, we develop
a novel scheme for solving nonparametric regression problems. That is the approximation of possibly low regular and noised functions from the knowledge of their approximate values given at some random points. Our proposed scheme is based  on the use of the pseudo-inverse  of a random projection matrix, combined with some specific properties of the Jacobi polynomials system, as well as some properties of positive definite random matrices. This scheme has the advantages to be stable, robust, accurate  and fairly  fast in terms of  execution time. In particular, we provide an $L_2$ as well as an $L_2-$risk errors of our proposed nonparametric regression estimator. Moreover and unlike most of the existing nonparametric regression estimators, no extra regularization step is required by our proposed estimator. Although, this estimator is initially designed to work with random sampling set of uni-variate i.i.d. random variables following a Beta distribution, we show that it is still works for a wide range of  sampling distribution laws. Moreover, we briefly describe how our estimator can be adapted in order to handle the multivariate case of random sampling sets.
In the second part of this work, we extend the random pseudo-inverse scheme technique to build a stable and accurate estimator for solving  linear functional regression (LFR) problems. A dyadic decomposition approach is used to construct this last stable estimator for the LFR problem. Alaso, we give an $L_2-$risk error of our proposed LFR estimator. Finally,  the performance of the two  proposed estimators are illustrated  by various numerical simulations. In particular, a real dataset is used to illustrate the performance of our nonparametric regression estimator. \\


\section{Introduction} 

In the first part of this work, we consider the following nonparametric regression model:
\begin{equation}\label{problem1}
Y_i = f(X_i) +\varepsilon_i =\widetilde f(X_i),\quad 1\leq i\leq n,
\end{equation}
where $(X_i)_{1\leq i \leq n}$ are random variables (or inputs) with distribution $\rho_X$ and the noise terms $(\varepsilon_i)_{1\leq i \leq n}$ are i.i.d. real-valued centered  random variables with variance $\sigma^2.$  For 
a given  complete metric space $\mathcal X,$  a measurable  subset $\mathcal Y$ of $\mathbb R$ and a training set $\{(X_i, Y_i), 1\leq i\leq n\}$ in $\mathcal X\times \mathcal Y,$ we are interested in the construction of an estimator $\widehat f$ such that for $x\in \mathcal X,$ $\widehat f(x)$ is a good estimate  for  the corresponding regression function $f.$ Recall  that if the  $(X_i,Y_i)$ are  drawn from a joint probability measure $\rho$ on $\mathcal X \times \mathcal Y,$ then the regression function $f$ is given by 
\begin{equation*}
f(x)= \int_{\mathcal Y} y d \rho(y|x)= \mathbb E(Y|X=x),\quad x\in \mathcal X,
\end{equation*}
 For simplicity, we will assume that the $X_i$ take values in  the interval $I=[-1, 1].$ Nonetheless, our  results are easily adapted to a more general compact interval of $\mathbb R.$  Also, we show how our nonparametric regression estimator can be adapted in order to handle the more general case of  random sampling points $X_i\in \mathbb R^d,\, d\geq 2.$\\
 
Note that the  kernel ridge  regression (KRR) scheme is a widely used scheme for solving the regression problem modeled by \eqref{problem1}.  For more details on the KRR based estimator, the reader is referred for example  to \cite{Smale1, Smale2}. 
Other estimators for nonparametric regression problem have been proposed in the literature. To cite but a few, F. Comte and her co-authors have proposed
estimators based on  projections associated with some classical  orthogonal polynomials \cite{Comte}. Also, in the literature, there exists another class of convolution kernels with variable bandwidths based robust  estimators for the  nonparametric regression problem, see for example \cite{Benhenni, Zhao}.\\
 
Our proposed estimator for solving problem \eqref{problem1} is based on combining  the Random Sampling Consensus (RANSAC) iterative algorithm with 
a random pseudoinverse based  scheme for nonparametric regression. That is  at each iteration of the RANSAC algorithm,  
 a stable nonparametric regression  estimator is constructed by using the  pseudoinverse of a random projection matrix, associated with 
Jacobi orthogonal polynomials. Note that the  RANSAC procedure, see for example \cite{Sharon} consists in randomly
selecting $m < n$ measurements from model \eqref{problem1}. From these $m$ training data points $(X_i, Y_i),$ construct 
an estimator and compute the mean squared prediction error estimation. This procedure is iterated till an estimator with high consensus is obtained. We should mention that in the literature,   pseudoinverse scheme has been successfully used in a wide range of applications from different scientific area, including machine learning, see \cite{Jung} and the references therein. The usual way of computing the pseudoinverse of an $n\times m$ real coefficients matrix is to combine the Moore-Penrose inverse and the SVD techniques.  It is well known that the time complexity of this last procedure is $O\big(\min(n m^2, n m^2)\big).$ Nonetheless in \cite{Jung}, the authors have proposed an accurate and  Fast PesudoInverse algorithm for sparse feature matrices.

For the sake of simplicity, the emphasis will be on the construction of a stable estimator for problem \eqref{problem1} by using a random training set 
$\{(X_i,Y_i),\,\, 1\leq i\leq n\}.$  More precisely,  for  real parameters ${\displaystyle \alpha, \beta \geq -\frac{1}{2},}$ we  consider 
the weight function $\omega_{\alpha,\beta}(x)= (1-x)^\alpha (1+x)^\beta$ and the associated weighted $L^2_\omega(I)-$space. This last space 
is  given by the set of real valued functions $f$ that are measurable on $I$ and satisfying
${\displaystyle \|f\|_\omega^2= \int_I |f(x)|^2 \omega_{\alpha,\beta}(x)\, dx  < +\infty.}$ 
For a positive integer  $n\in \mathbb N,$ let $\{ X_i,\, 1\leq i\leq n\}$ be a  sampling set of i.i.d random variables following the beta distribution $B(\alpha+1,\beta+1)$ over $I.$ Note that this last condition is not mandatory and it is used in order to simplify the different theoretical results related 
to our nonparametric regression estimator. In fact, we show how this estimator can be adapted in order to handle random sampling set drawn from a more general probability distribution law or even from an unknown probability distribution.\\

Our  nonparametric regression  estimator is briefly described as follows.
We use the notation $\pmb u'$ for the transpose of $\pmb u,$ then for a positive integer $N\in \mathbb N$ with $N+1\leq n,$ 
$\widehat f_{n,N}$ is given in terms of the first $N+1$ normalized Jacobi polynomials $\wJ_{k}(x), \ k=0,\ldots, N.$ More precisely, we have 
\begin{equation}\label{Estimation_Jacobi1}
\widehat f_{n,N}(x)=\sum_{k=1}^{N+1} \widehat c_k \wJ_{k-1}(x),\quad x\in I.
\end{equation}
Here, the expansion coefficients vector $\widehat{\pmb C}= (\widehat c_1,\ldots, \widehat c_{N+1})'$ is given by the over-determined system 
\begin{equation}\label{B_N}
B_N \cdot  \widehat{\pmb C} = \Big[\frac{1}{\sqrt{n}} \wJ_{k-1}(X_j) \Big]_{j,k} \cdot \widehat{\pmb C} =\mathbf Z_n,
\end{equation}
where  $\mathbf Z_n = \frac{1}{\sqrt{n}} (Y_1,\ldots,Y_n)'.$  We refer to the $n\times (N+1)$ matrix  $B_N$ as the random projection matrix. One of the main result of this work is to prove that with high probability, 
the $N+1$ dimensional positive definite  matrix $A_N= B_N' B_N$ is invertible and well conditioned. More precisely, we show that with high probability, the $2-$norm condition number of $A_N,$
denoted by $\kappa_2(A_N)$ is bounded by a fairly small constant. This ensures the stability of the estimator  $\widehat f_{n,N}.$ Consequently, the expansion coefficients vector $\widehat{\pmb C}$ is simply given by 
\begin{equation}\label{coefficients}
\widehat{\pmb C}= \big( B_N' B_N\big)^{-1} B_N' \mathbf Z_n.
\end{equation}
Also, we provide an upper bound for the regression error $\| f-\widehat f_{n,N} \|_\omega,$ that holds with high probability.   In particular, we provide more precise regression  errors in the special cases where  $f$ has a  Lipschitzian   $p-$th order derivative or $f$ is a  $c-$bandlimited function, for some $c>0.$ Moreover, we provide a bound for the 
$L^2-$risk error  ${\displaystyle \mathbb E\Big[\| f-\widehat f_{n,N} \|^2_\omega\Big]}.$

In the second part of this work, we push forward the random pseudo-inverse based scheme  for solving a linear functional regression (LFR) problem, that was introduced in \cite{BDK}. This is done by combining the random  pseudo-inverse scheme with a dyadic decomposition technique. This combination yields a stable and a  highly accurate  estimator for solving the LFR problem. Recall that for a compact interval $J,$ the LFR model is given as follows
\begin{equation}\label{problem2}
  Y_i =\int_J X_i(s)\,  \beta_0(s)\, ds +\varepsilon_i,\quad i=1,\ldots,n.
 \end{equation}
Here, the $X_i(\cdot)\in L^2(J)$ are random functional predictors, the $\varepsilon_i$ are i.i.d. centered white noise independent  of the $X_i(\cdot)$ and 
$\beta_0(\cdot)\in L^2(J)$ is the unknown slope function to be recovered. As usual, we assume that the  $X_i(\cdot)$ are given by their series expansions with respect to an orthonormal set of $L^2(J).$  That is for an orthonormal family  ${\displaystyle \{\varphi_k(\cdot),\, k \geq 1\}}$ of $L^2(J),$
 we have ${\displaystyle X_i(s)= \sum_{k\geq 1} \xi_k Z_{i,k} \varphi_k(s),}$ where
 the   $Z_{i,k}$ are i.i.d. centered  random variables with variance $\sigma_Z^2$ and $(\xi_k)_{k\geq 1}$ is a deterministic  sequence of $\mathbb R\setminus  
 \{0\}.$  Our proposed stable random pseudo-inverse based estimator  $\widehat \beta_{n,N}(\cdot)$ of ${\displaystyle \beta_0(\cdot)=\sum_{k\geq 1} c_k \varphi_k(\cdot)}$ is briefly described as follows.
For an integer $N\geq 2,$ we let $K_N=\Big[\frac{\log N}{\log 2}\Big]+1$ and we consider the dyadic decomposition of the 
set $[[1,N]]=\{1,\ldots,N\}$ into $K_N$ subsets   $ I_k=[[2^{k-1}+1,\min(2^k,N)]],$ $1\leq k\leq K_N.$ 
On each subset 
	$I_k,$ we consider the random matrices  $G_k, F_k$ given by 
	$$G_k=\displaystyle\left[\frac{1}{n} \sum_{i=1}^n \xi_j\xi_l Z_{i,j}Z_{i,l} \right]_{j,l\in I_k}= F_k' F_k,\quad F_k=\frac{1}{\sqrt n}\big[\xi_jZ_{i,j}\big]_{\substack{ \scriptstyle 1\leq i\leq N\\ \scriptstyle   j\in I_k}}.$$
	Under the previous notation, our proposed estimator for the LFR problem \eqref{problem2} is given by 
$$\widehat \beta_{n,N}(s)=\sum_{k=1}^{K_N} \sum_{j\in I_k} \widehat c_{n,j} \varphi_j(s),\quad \widehat{\pmb c}_{n,k} =[\widehat c_{n,j}]_{j\in I_k}'=G_k^{-1}\Big(F_k'\frac{1}{\sqrt{n}}[Y_i^k]'_{1\leq i\leq n}\Big),\quad 1\leq k \leq K_N.$$
Note that the  $\pmb Y^k = [Y_i^k]_{1\leq i\leq n}$ are obtained from the output vector $\pmb Y=[Y_i]_{1\leq i\leq n}$ by substituting in \eqref{problem2},  $X_i(s)$ with its projection ${X_i^k(s)=\sum_{j\in I_k} \xi_j Z_{i,j} \varphi_j(s)}.$
In particular, we provide an upper bound for the mean squared estimation error $ \|\widehat{\beta}_{n,N}(\cdot)-\beta_0(\cdot)\|_{2}^2,$ that holds with high probability. This upper error bound involves the quantity ${\displaystyle \sum_{k=1}^{K_N}\kappa_2(G_k).}$  In practice, the deterministic sequence $\pmb \xi =(\xi_j)_j$ satisfies a decay condition of the type $\xi_j \asymp j^{-s}$ for some $s>0.$ Under this last assumption, it is easy to check that our combined   dyadic decomposition and random pseudo-inverse scheme is stable in the sense that ${\displaystyle \sum_{k=1}^{K_N}\kappa_2(G_k) \lesssim 2^s \log N,}$ which is a relatively  small cumulative condition number. This ensures the stability of our proposed LFR estimator $\widehat{\beta}_{n,N}(\cdot).$ Also, we give an estimate for the $L^2-$risk error  $\mathbb E\Big[ \| \beta_0(\cdot)-\widetilde \beta_{N,M}(\cdot) \|^2_2\Big],$ where $\widetilde \beta_{N,M}(\cdot)$ is a truncated version of the estimator $\widehat{\beta}_{n,N}(\cdot).$\\

This work is organized as follows. In section 2, we give some mathematical preliminaries that are frequently used in this work. In section 3, we study our random pseudo-inverse based estimator for the stable and robust  approximation of the regression function $f$ given by \eqref{problem1}. In section 4, we give an error analysis of our combined dyadic decomposition, random pseudo-inverse based estimator for the LFR problem \eqref{problem2}. In section 5, we give some numerical simulations that illustrate the different results of this work. In particular, we illustrate the performance of our nonparametric regression by applying it on a real dataset given by the daily Covid-19 cases of some world's countries.

\section{Mathematical preliminaries}

We first   recall that for two real numbers $\alpha, \beta >-1,$ the classical Jacobi polynomials $\J_k$ are defined for $x\in [-1,1],$ by the following Rodrigues formula,
\begin{equation}
\label{Jacobi_Rodrigues}
\J_k (x)= \frac{(-1)^k}{2^k k!} \frac{1}{\omega_{\alpha,\beta}(x)} \frac{\mbox{d}^k}{\mbox{d} x^k} \Big( \omega_{\alpha,\beta}(x)(1-x^2)^k\Big),\quad
\omega_{\alpha,\beta}(x)= (1-x)^\alpha (1+x)^\beta,
\end{equation}
with
$$\J_k(1) = { k+\max(\alpha,\beta) \choose{k}}  =\frac{\Gamma(k+\max(\alpha,\beta)+1)}{k!\, \Gamma(\max(\alpha,\beta)+1)}.$$
Here, $\Gamma(\cdot)$ denotes the usual Gamma function. 
In the sequel, we let $\wJ_k$ denote the normalized Jacobi polynomial of degree $k,$ given by 
\begin{equation}\label{JacobiP}
\wJ_{k}(x)= \frac{1}{\sqrt{h_k^{\alpha,\beta}}}\J_k(x),\quad h_k^{\alpha,\beta}=\frac{2^{\alpha+\beta+1}\Gamma(k+\alpha+1)\Gamma(k+\beta+1)}{k!(2k+\alpha+\beta+1)\Gamma(k+\alpha+\beta+1)}.
\end{equation}
In this case, we have
\begin{equation}\label{Normalisation1}
\| \wJ_k \|^2_\omega=\int_{-1}^1 (\wJ_k(y))^2 \omega_{\alpha,\beta}(y)\, dy =1.
\end{equation}
The following lemma provides us with some useful upper bounds for the normalized Jacobi polynomials $\wJ_k,$ for $k\geq 2.$

\begin{lemma} For $\alpha, \beta \geq -\frac{1}{2},$   let $\mu=\max(\alpha,\beta),$ $c_{\alpha,\beta}=\frac{\alpha+\beta+1}{2}$ and let 
\begin{equation}
\label{Eq2.1}
 \eta_{\alpha,\beta}=\frac{\exp\big(\frac{2\max(\mu,0)}{12}+\frac{\max(\mu^2+\alpha \beta,0)}{8}\big)}{ 2^{(\alpha+\beta)/2}\Gamma(\mu+1)}
\end{equation}
Then, we have 
\begin{equation}
\label{Ineq2.1}
\max_{x\in [-1,1]} |\wJ_k (x)| \leq \eta_{\alpha,\beta} \gamma_k^{\alpha,\beta},\qquad \gamma_k^{\alpha,\beta}=   k^{\mu} \sqrt{k+c_{\alpha,\beta}},\quad \forall\, k\geq 2.
\end{equation}

\end{lemma}  

\noindent{\bf Proof:} It is well known that under the previous condition on $\alpha, \beta,$ we have 
\begin{equation}
\label{Ineq2.2}
\max_{x\in [-1,1]} |\J_k(x)|= \max\big(|\J_k(1)|,|\J_k(-1)|\big) = \frac{\Gamma(k+\mu+1)}{k!\, \Gamma(\mu+1)},\quad \mu=\max(\alpha,\beta).
\end{equation}
On the other hand, from \cite{Wang}, for an integer $n\geq 2$ and for two real numbers $a,b>-1,$ we have
\begin{equation}
\label{Ineq2.3}
\frac{\Gamma(n+a)}{\Gamma(n+b)} \leq C_n^{a,b} n^{a-b},\quad C_n^{a,b}= \exp\Big(\frac{a-b}{2(n+b-1)}+\frac{1}{12(n+a-1)}+\frac{(a-1)(b-1)}{n}\Big).
\end{equation}
By combining \eqref{JacobiP}, \eqref{Ineq2.2} and \eqref{Ineq2.3}, one  gets
$$\max_{x\in [-1,1]} |\wJ_k (x)|\leq \frac{\sqrt{C_k^{\mu+1,1}C_k^{\alpha+\beta+1,\alpha+\beta-\mu+1}}}{\Gamma(\mu+1) 2^{(\alpha+\beta)/2}}\, \,  k^\mu\,\, \sqrt{k+\frac{(\alpha+\beta+1)}{2}}.$$
But from \eqref{Ineq2.3}, it can be easily checked that for $k\geq 2,$ we have 
$$C_k^{\mu+1,1}C_k^{\alpha+\beta+1,\alpha+\beta-\mu+1}\leq \eta_{\alpha,\beta}.$$
By combining the previous two inequalities, we get
\begin{equation}\label{Ineq2.4}
\max_{x\in [-1,1]} |\wJ_k (x)| \leq \eta_{\alpha,\beta}  \gamma_k^{\alpha,\beta}= \eta_{\alpha,\beta} k^\mu\,\, \sqrt{k+\frac{(\alpha+\beta+1)}{2}},\quad k\geq 2.
\end{equation}

Next, we recall from \cite{Johnson}, the following  Gershgorin circle theorem that is used for the estimation of $2-$norm condition numbers of various random
matrices used in this work.

\begin{theorem}[Gershgorin circle theorem] Let $A= [a_{ij}]$ be a complex $n\times n$  matrix. For $1\leq i\leq n,$ let 
${\displaystyle R_i = \sum_{j\neq i} |a_{ij}|.}$ Then, every eigenvalue of $A$ lies within at least one of the discs $D(a_{ii}, R_i).$
\end{theorem}

The following McDiarmid's concentration inequality is also used in the spectral analysis of the random matrices of the next section.

\begin{theorem}[McDiarmid's inequality] Let $X_1,\ldots,X_n$ be independent  random variables with respective ranges in metric spaces 
$\pmb\chi_1,\ldots,\pmb\chi_n.$ Let $f:\pmb \chi_1\times\cdots\times\pmb \chi_n\rightarrow \mathbb C$ be a function satisfying the bounded differences property 
$$\Big| f(x_1,\ldots,x_{i-1},x_i,x_{i+1},\ldots,x_n)-f(x_1,\ldots,x_{i-1},x'_i,x_{i+1},\ldots,x_n)\Big|\leq c_i,\quad 1\leq i\leq n.$$
Then, for any $\epsilon >0,$ we have 
$$\mathbb P \Big( f(X_1,\ldots,X_n)-\mathbb E\big[ f(X_1,\ldots,X_n)\big]\leq \epsilon \Big)\geq 1-\exp\left(\frac{2\epsilon^2}{\sum_{i=1}^n c_i^2}\right).$$
\end{theorem}

\section{ Random projection on Jacobi polynomials system and nonparametric regression estimator}

The main part of this section is devoted to the study of  the stability and error analysis of our nonparametric regression estimator in the framework
of a sampling set drawn from a Beta distribution. Then, at the end of this section, we show how our estimator can be adapted in order to handle random sampling sets drawn from more general distribution laws or even from an unknown distribution law. \\

We first consider  two real numbers $\alpha, \beta \geq -\frac{1}{2}$ and  two positive integers 
$n\geq N+1$, then we   let $\{X_i,\, 1\leq i\leq n\}$ be a sampling set of i.i.d random variables following the beta $B(\alpha+1,\alpha+1)$ distribution over the interval $I.$ Then, under the assumption that the regression function $f$ (solution of the nonparametric regression problem \eqref{problem1}) is well approximated by its projection $f_N =\pi_N f$ over the subspace  $\mathcal H_N=\mbox{Span}\{ \wJ_{k},\, 0\leq k\leq N\},$
our random pseudo-inverse based estimator $\widehat f_{n,N}$ of $f\in L^2(I)$ is given by 
\begin{equation}\label{Estimation_Jacobi}
\widehat f_{n,N}(x)=\sum_{k=1}^{N+1} \widehat c_k \wJ_{k-1}(x),\quad x\in I.
\end{equation}
We require that  this estimator is stable and robust. This estimator is given as the Least square solution of the over-determined system    
\begin{equation}\label{Estimation_Jacobi2}
\sum_{k=1}^{N+1} c_k \wJ_{k-1}(X_j)= Y_j,\quad j=1,\ldots,n.
\end{equation}
Note that if the regression function $f\in L^2_\omega(I),$ then  from the uniqueness of its expansion coefficients with respect to the orthonormal 
basis $\{ \wJ_k(\cdot),\, k\geq 0\},$ for any integer $N\in \mathbb N,$ its  first $N$ coefficients $c_i(f)$  do not depend on the specific values
of $y_j = f(X_j),$ at a specific random sampling set $\{ X_j, j=1,\ldots,n\}.$ These coefficients are uniquely given by
$$ c_k(f) = < \widetilde f, \wJ_{k-1}> = \int_I \widetilde f(x) \wJ_{k-1}(x) \omega_{\alpha,\beta}(x)\, dx,\quad 1\leq k\leq N+1.$$
Note that by re-scaling by a factor ${ \frac{1}{\sqrt{n}},}$  the over-determined system \eqref{Estimation_Jacobi2} is written as
\begin{equation}\label{system2}
B_N \cdot  \pmb C = \pmb Y,\quad B_N= \Big[\frac{1}{\sqrt{n}} \wJ_{k-1}(X_j) \Big]_{1\leq j\leq n,1\leq k\leq N+1},\quad \pmb Y= [Y_j]'_{1\leq j\leq n}.
\end{equation}
We show that under some conditions on $f$ as well as on the integers $n, N,$ the least square solution of the previous system is a  good  approximation
of the set of the first $N+1$ expansion coefficients of $f$. For this purpose, we define an $(N+1)\times (N+1)$ random Jacobi projection matrix $A_N$ by
\begin{equation}\label{Ineq1.2}
A_N = B_N' B_N = \left[\frac{1}{n} \sum_{j=1}^n \wJ_{k-1}(X_j)\wJ_{l-1}(X_j) \right]_{1\leq k,l\leq N+1} = \sum_{j=1}^n \mathbf D_j,
\end{equation}
Note that the random matrix $\mathbf D_j$ is also positive semidefinite. In fact, $\mathbf D_j = b_j' b_j,$ where $b_j$ is the $1\times (N+1)$ matrix given by 
${\displaystyle b_j= \frac{1}{\sqrt{n}}\Big[ \wJ_{0}(X_j) \cdots  \wJ_{N}(X_j)\Big]}.$ Consequently, for any $X\in \mathbb R^{N+1},$ we have 
$$ X' \cdot \mathbf D_j X = X' b_j' \cdot  b_j X  \geq 0.$$ That is $\mathbf D_j$ is a positive semi-definite matrix.
An important property of the matrix $A_N$ is that under some conditions on the parameters $\alpha, \beta, n, N,$ this matrix is well conditioned.
To prove this first main result, we need the following proposition that provides us with 
an upper and a lower bound for $\lambda_{\max}(A_N)$ and $\lambda_{\min}(A_N),$ respectively.

\begin{proposition}\label{proposition1} Under the notation of Lemma 1, for any two positive integers $n> N\geq 2,$ satisfying 
\begin{equation}\label{Condition1}
0.63\,  n >  m_{\alpha,\beta}^2 \cdot  (N+1)^{2\mu+2}\, \log(N+1),
\end{equation}
where ${\displaystyle  m_{\alpha,\beta}^2 = \frac{1+\frac{1}{2}\sqrt{\frac{c_{\alpha,\beta}}{2}}}{\mu+3/2} \eta^2_{\alpha,\beta}}$ and $\eta_{\alpha,\beta}$ is as given by \eqref{Eq2.1}. 
Then, we have 
\begin{equation}\label{Expectation1}
\mathbb E\big(\lambda_{\max}(A_N)\big)\leq 1.72 + m_{\alpha,\beta}^2\frac{(N+1)^{2\mu+2}\, \log(N+1)}{n} 
\end{equation}
and
\begin{equation}\label{EExpectation1}
\mathbb E\big(\lambda_{\min}(A_N)\big)\geq 0.63-m_{\alpha,\beta}^2 \frac{(N+1)^{2\mu+2}\, \log(N+1)}{n}.
\end{equation}
\end{proposition}

\noindent
{\bf Proof:} 
We first note that from \eqref{JacobiP}, we have  
\begin{eqnarray*}
\max_{x\in [-1,1]} |\wJ_0 (x)|&=&\big(2^{\alpha+\beta+1} B(\alpha+1,\beta+1)\big)^{-1/2}\\
\max_{x\in [-1,1]} |\wJ_1 (x)| &= &(1+\mu) \sqrt{\frac{3+\alpha+\beta}{(\alpha+1)(\beta+1)}}\big(2^{\alpha+\beta+1} B(\alpha+1,\beta+1)\big)^{-1/2}.
\end{eqnarray*}
On the other hand, if $\gamma_k^{\alpha,\beta}$ is as given by \eqref{Ineq2.1}, then it  is easy to check that for $\mu\geq -\frac{1}{2}$ and for $\alpha, \beta \geq -\frac{1}{2},$ the function 
$x\rightarrow x^\mu \sqrt{x+c_{\alpha,\beta}}$ is increasing on $[2, +\infty).$ Consequently, we have 
\begin{eqnarray*}
\sum_{k=2}^N \gamma_k^{\alpha,\beta}&\leq & \int_2^{N+1} x^\mu \sqrt{x+c_{\alpha,\beta}}\, dx \leq \int_2^{N+1} x^{\mu+\frac{1}{2}} \Big(1+\frac{1}{2}\sqrt{\frac{c_{\alpha,\beta}}{2}}\Big)\, dx\\
& \leq & \frac{\Big(1+\frac{1}{2}\sqrt{\frac{c_{\alpha,\beta}}{2}}\Big)}{\mu+3/2} (N+1)^{\mu+3/2}.
\end{eqnarray*}
Hence, by using \eqref{Ineq2.1} and  straightforward computation, one gets 
$$ \max_{x\in [-1,1]} |\wJ_0 (x)|+\max_{x\in [-1,1]} |\wJ_1 (x)|+\cdots+\max_{x\in [-1,1]} |\wJ_N (x)| \leq \eta_{\alpha,\beta} \int_N^{N+1} x^\mu \sqrt{x+c_{\alpha,\beta}}\, dx.$$
The previous two inequalities give us
\begin{eqnarray}\label{Ineq2.5}
\lefteqn{\max_{x\in [-1,1]} |\wJ_N (x)| \sum_{k=1}^{N+1} \Big(\max_{x\in [-1,1]} |\wJ_{k-1} (x)|\Big)\leq \eta^2_{\alpha,\beta} N^\mu \sqrt{N+c_{\alpha,\beta}}  \frac{\Big(1+\frac{1}{2}\sqrt{\frac{c_{\alpha,\beta}}{2}}\Big)}{\mu+3/2} (N+1)^{\mu+3/2}}\nonumber \\
& &\qquad\qquad\quad\qquad\qquad\qquad\leq  \frac{\Big(1+\frac{1}{2}\sqrt{\frac{c_{\alpha,\beta}}{2}}\Big)}{\mu+3/2} \eta_{\alpha,\beta}^2 (N+1)^{2\mu+2}=m_{\alpha,\beta}^2 (N+1)^{2\mu+2}.
\end{eqnarray}
Next, let $\mathbf D_j$ be the positive semi-definite matrix given by \eqref{Ineq1.2}. Since the random samples $X_i$ follow the Beta-distribution on $[-1,1]$
with parameters $(\alpha+1,\beta+1),$ then we have $\mathbb E \left(\wJ_{k-1}(X_j)\wJ_{k-1}(X_j) \right)= \delta_{kl}.$ That is 
${\displaystyle \mathbb E\Big(\sum_{j=1}^n \mathbf D_j\Big)= I_{N+1},}$ where $I_{N+1}$ is the identity matrix of dimension $N+1.$ Hence, we have 
$$ \Lambda_{\min} = \lambda_{\min}\Big(\mathbb E\big(\sum_{j=1}^n \mathbf D_j\big)\Big)=1,\quad \Lambda_{\max} = \lambda_{\max}\Big(\mathbb E \big(\sum_{j=1}^n \mathbf D_j\big)\Big)=1.$$
Also, from  Gershgorin circle theorem, we have 
\begin{equation}
\label{Ineq3.1}
\lambda_{\max}(\mathbf D_j) \leq \frac{1}{n} \max_{x\in [-1,1]} |\wJ_N (x)| \sum_{k=1}^{N+1} \Big(\max_{x\in [-1,1]} |\wJ_{k-1} (x)|\Big)\leq \frac{1}{n} m_{\alpha,\beta}^2 (N+1)^{2\mu+2}=\frac{L_N(\alpha,\beta)}{n}.
\end{equation}
Moreover, from \cite{Tropp}, we have 
\begin{equation}\label{Ineq2.6}
\mathbb E\Big(\lambda_{\min} (\sum_{j=1}^n \mathbf D_j\big)\Big)\geq 0.63 \Lambda_{\min} - \frac{L_N(\alpha,\beta)}{n} \log(N+1),\quad
\mathbb E\Big(\lambda_{\max} \big(\sum_{j=1}^n \mathbf D_j\big)\Big)\leq 1.72 \Lambda_{\max} + \frac{L_N(\alpha,\beta)}{n} \log(N+1).
\end{equation}
This concludes the proof of the proposition.\\

\noindent
The following theorem provides us with an upper bound for the actual condition number $\kappa_2(A_N).$
\begin{theorem}\label{thm1}
Under the previous notation and the same conditions as proposition \ref{proposition1}, for any $0<\delta<1,$ we have  
\begin{equation} \label{deviation1}
\kappa_2 (A_N)\leq \frac{1.72+m_{\alpha,\beta}^2 (N+1)^{2\mu+2}\, \Big(\frac{\log(N+1)}{n}+\sqrt{\frac{2}{n}\log\big(\frac{2}{\delta}\big)}\Big)}{0.63-m_{\alpha,\beta}^2 (N+1)^{2\mu+2}\, \Big(\frac{\log(N+1)}{n}+\sqrt{\frac{2}{n}\log\big(\frac{2}{\delta}\big)}\Big)}
\end{equation}
with probability at least $(1-\delta)^2.$  
\end{theorem}

\noindent
{\bf Proof:}
Given  $x=(x_1,\cdots,x_n)\in [-1,1]^n$, we consider the $(N+1)\times(N+1)$ matrix with entries $\frac{1}{n} \sum_{j=1}^n \wJ_{k-1}(x_j)\wJ_{l-1}(x_j)$. We use the notation $A_N(x)$ for such a matrix.
We denote by $\lambda_k(A_N(x))$, $1\leq k \leq N+1,$ the eigenvalues of $A_N(x)$ arranged in decreasing order.
For $1\leq k \leq N+1,$ we use McDiarmid's concentration inequality for the $n-$variate mapping 
$$x \mapsto \phi_k(x)=\lambda_k(A_N(x)) $$
We  prove that the previous mapping $\phi_k$ satisfies the bounded differences  assumption.
That is, when only one of the $n$ coordinates differs between $x$ and $x'$, then
\begin{equation}
\label{Diarmid1}
|\phi_k(x)-\phi_k(x')| \leq \frac{2}{n} m_{\alpha,\beta}^2 (N+1)^{2\mu+2}=2\frac{L_N(\alpha,\beta)}{n}.
\end{equation}
Let us take this for granted for now and we will prove it later on.  It follows from McDiarmid's inequality that 
\begin{equation} \label{Ineq2.8}
\mathbb P\Big ( \phi_k(X_1,\cdots,X_n)-\mathbb E\phi_k(X_1,\cdots,X_n) \geq \varepsilon\Big) \leq \exp\frac{-n\varepsilon^2}{2m_{\alpha,\beta}^4 (N+1)^{4\mu+4}}=\exp\frac{-n\varepsilon^2}{2L_N^2(\alpha,\beta)}
\end{equation}
and
\begin{equation} \label{Ineq2.9}
\mathbb P\big ( \phi_k(X_1,\cdots,X_n)-\mathbb E\phi_k(X_1,\cdots,X_n) \leq -\varepsilon\big) \leq \exp\frac{-n\varepsilon^2}{2L_N^2(\alpha,\beta)}.
\end{equation}
For the value of $k=1$, \eqref{Ineq2.8} gives us
\begin{equation} \label{Ineq2.10}
\mathbb P\big ( \lambda_{\max}(A_N))-\mathbb E(\lambda_{\max}(A_N))   \leq  \varepsilon\big) \geq 1- \exp\frac{-n\varepsilon^2}{2L_N^2(\alpha,\beta)}
\end{equation}
Moreover, for the value of $k=N+1$, \eqref{Ineq2.9} gives  us
\begin{equation} \label{Ineq2.11}
\mathbb P\big ( \lambda_{\min}(A_N)-\mathbb E(\lambda_{\min}(A_N)) \geq -\varepsilon\big) \geq 1- \exp\frac{-n\varepsilon^2}{2L_N^2(\alpha,\beta)}
\end{equation}
By combining \eqref{Ineq2.6} and \eqref{Ineq2.10}, one gets 
\begin{equation} \label{Ineq2.12}
\lambda_{\max}(A_N) \leq  1.72+\frac{L_N(\alpha,\beta)}{n}\log(N+1)+\varepsilon
\end{equation}
with probability at least $1-\exp\frac{-n\varepsilon^2}{2L_N^2(\alpha,\beta)}$. Here, $L_N(\alpha,\beta)$ is given by \eqref{Ineq3.1}. Also, 
by combining \eqref{Ineq2.6} and \eqref{Ineq2.11}, one gets
 \begin{equation} \label{Ineq2.13}
\lambda_{\min}(A_N)\geq \big(0.63- \frac{L_N(\alpha,\beta)}{n} \log(N+1) \big)-\varepsilon  
\end{equation}
with probability at least $1-\exp\frac{-n\varepsilon^2}{2L_N^2(\alpha,\beta)}$.
Hence, by combining \eqref{Ineq2.12} and \eqref{Ineq2.13} , one gets 
$$\kappa_2 (A_N)\leq \frac{1.72+\frac{L_N(\alpha,\beta)}{n}\log(N+1))+\varepsilon}{0.63- \frac{L_N(\alpha,\beta)}{n} \log(N+1)-\varepsilon }$$
with probability at least ${\displaystyle \Big(1-\exp\frac{-n\varepsilon^2}{2L_N^2(\alpha,\beta)}\Big)^2}$. Finally to get \eqref{deviation1}, it suffices to let 
${\displaystyle \delta = \exp\left(-\frac{n\varepsilon^2}{2 L_N^2(\alpha,\beta)}\right)}.$\\

\noindent
Let us now prove the bounded differences  condition \eqref{Diarmid1}. We have  
\begin{eqnarray*}
|\phi_k(x)-\phi_k(x')|&=&|\phi_k(x_1,\ldots,x_{i-1},x_i,x_{i+1}\cdots,x_n)-\phi_k(x_1,\ldots,x_{i-1},x_i',x_{i+1}\cdots,x_n)|\\
&=&\lambda_k(A_N(x))-\lambda_k(A_N(x'))
\end{eqnarray*}
Let $E$ be the $(N+1)\times(N+1)$ matrix with entries $$E_{k,l}=\frac{1}{n}\left(\wJ_{k-1}(x_i)\wJ_{l-1}(x_i)-\wJ_{k-1}(x'_i)\wJ_{l-1}(x'_i)\right),$$
so that $A_N(x)=A_N(x')+E$. From Weyl's perturbation theorem of the spectrum of a perturbed Hermitian matrix, see for example \cite{Braun}, we have
$$\lambda_k(A_N(x'))+\lambda_{N+1}(E)\leq \lambda_k(A_N(x')+E)\leq \lambda_k(A_N(x'))+\lambda_1(E), \quad k=1,\cdots,N+1.$$
That is  $$|\lambda_k(A_N(x'))-\lambda_k(A_N(x))|\leq ||E||.$$
Moreover, from Gershgorin circle theorem, we have, for $1\leq j\leq N+1$
$$|\lambda_j(E)-E_{j,j}|\leq \sum_{\substack{1\leq p \leq N+1\\p\neq j}}|E_{j,p}|.$$
Hence, one gets 
\begin{eqnarray*}
|\lambda_j(E)|&\leq& |E_{j,j}|+\sum_{\substack{1\leq p \leq N+1\\p\neq j}}|E_{j,p}|=\sum_{1\leq p\leq N+1}|E_{j,p}|\\
&\leq&\frac{1}{n} \sum_{1\leq p\leq N+1} \big|\wJ_{j-1}(X_i)\wJ_{p-1}(X_i)|+|\wJ_{j-1}(X'_i)\wJ_{p-1}(X'_i)\big|\\
&\leq & \frac{2}{n} \max_{x\in [-1,1]} |\wJ_N (x)|\sum_{k=1}^{N+1} \Big(\max_{x\in [-1,1]} |\wJ_k (x)|\Big)\leq \frac{2}{n} m_{\alpha,\beta}^2 (N+1)^{2\mu+2}=2\frac{L_N(\alpha,\beta)}{n}.
\end{eqnarray*}
This last inequality follows from \eqref{Ineq2.5}.
Since $E$ is Hermitian, then  $||E||=\displaystyle\max_{\substack{1\leq j\leq N+1}}|\lambda_j|$ which concludes the proof.

\begin{remark}\label{rem2} For the special convenient values of $\alpha=\beta=-\frac{1}{2},$ the estimate of
$\kappa_2 (A_N)$ given by the previous theorem  can be further improved. In fact, in this case,
the normalized Jacobi polynomials are reduced to the normalized Chebyshev polynomials $\widetilde T_n$ defined on $[-1,1]$ by
$$\widetilde T_n(x) = \frac{2}{\pi} T_n(x),\quad T_n(cos\theta)= \cos(n\theta),\quad \theta \in [0,\pi].$$
It is easy to see that in this case, the optimal value of ${\displaystyle m_{-1/2,-1/2}=  \frac{2}{\pi}.}$ Consequently, with high probability, we have 
 \begin{equation}\label{Expectation1-2}
 \kappa_2 (A_N) \lesssim \frac{1.72 + \frac{2}{\pi}\frac{(N+1) \, \log(N+1)}{n}}{0.63- \frac{2}{\pi}\frac{(N+1)\, \log(N+1)}{n}}.
 \end{equation}
\end{remark}

\begin{remark}\label{rem2_2} Note that the upper bound of $\kappa_2(A_N)$ given by 
\eqref{deviation1} of the previous theorem is the essential condition for the validity 
of our nonparametric regression estimator. Numerical results given in the last section of this work,  indicate that the quantity $\kappa_2(A_N)$ remains bounded by a fairly convenient positive constant when the i.i.d. random sampling points follow a more general probability law. That is our proposed nonparametric regression estimator can be used in the framework of a more general random sampling set $\{X_i,\, 1\leq i\leq n\}.$
\end{remark}

Next, to study the regression estimation  error in the $\|\cdot\|_\omega-$norm, we need the following technical lemma. 

\begin{lemma}\label{Lemme1} Let $g\in L^2_\omega(I)$ be a  bounded function. Under the previous notations, for any $\alpha\geq -\frac{1}{2}$ and any $0<\delta <1,$
we have with probability at least $1-\delta,$
\begin{equation}\label{L2norm}
\left|\frac{\gamma_{\alpha,\beta}}{n} \sum_{j=1}^n \big(g(X_j)\big)^2-\| g\|_\omega^2\right|\leq \frac{\gamma_{\alpha,\beta}}{\sqrt{n}}  \sqrt{\log\big(\frac{2}{\delta}\big)} \|g\|_\infty^2.
\end{equation}
where $\gamma_{\alpha,\beta}=B(\alpha+1,\beta+1)2^{\alpha+\beta+1}.$
\end{lemma}

\noindent
{\bf Proof:} It suffices to use McDiarmid's inequality with the real valued $n-$variate function,
$$\phi: I^n \rightarrow \mathbb R,\,\, (x_1,\ldots,x_n)\rightarrow \frac{\gamma_{\alpha,\beta}}{n} \sum_{i=1}^n \big(g(x_i)\big)^2.$$
Then, for any $x_1,\ldots,x_{i-1}, x'_i,x_{i+1},\ldots,x_n\in I,$ we have
\begin{eqnarray*}
\lefteqn{|\phi(x_1,\ldots,x_{i-1},x_i,x_{i+1},\ldots,x_n)-\phi(x_1,\ldots,x_{i-1},x'_i,x_{i+1},\ldots,x_n)|  }\\
&&\qquad\qquad\qquad\qquad\qquad\qquad\qquad\qquad =\frac{\gamma_{\alpha,\beta}}{n}  \Big|\big(g(x_i)\big)^2-\big(g(x'_i)\big)^2\Big | \leq 
\frac{\gamma_{\alpha,\beta}}{n} \|g\|_\infty^2 =d_i.
\end{eqnarray*}
Consequently, we have ${\displaystyle \sum_{i=1}^n d_i^2 =\frac{(\gamma_{\alpha,\beta})^2}{n} \|g\|_\infty^4.}$ On the other hand, we have 
$$ \mathbb E \left(\frac{\gamma_{\alpha,\beta}}{n}\sum_{i=1}^n \big(g(X_i)\big)^2 \right)= \mathbb E \Big( \gamma_{\alpha,\beta}\big(g(X)\big)^2\Big)=\int_I |g(x)|^2\, \omega_{\alpha,\beta}(x)\, dx=\| g\|_\omega^2.$$

The following theorem provides us with an upper bound in the $\|\cdot\|_\omega-$norm for the regression error in terms 
of $\|\pi_N f\|_\omega$ and $\|\pi_N f\|_\infty.$ Here, $\pi_N f $ is the orthogonal  projection of the regression function over
the subspace $\mathcal H_N=\mbox{Span}\{\wJ_k(\cdot),\, 0\leq k\leq N\}.$ 

\begin{theorem} Under the previous notations and the hypotheses of Theorem 3, let ${\displaystyle \pmb \eta_n=\max_{1\leq i\leq n} |\varepsilon_i|}$, then  for any $0<\delta <1,$ we have with probability at least 
$(1-\delta)^2,$
\begin{equation}\label{Error1}
\| f-\widehat f_{n,N} \|_\omega \leq \|f-\pi_N f\|_\omega+ \sqrt{2 \kappa_2(A_N)} \frac{{\Big(\frac{1}{n}\log\big(\frac{2}{\delta}\big)\Big)}^\frac{1}{4}\|f-\pi_N f\|_\infty +\frac{\| f-\pi_N f\|_\omega}{\sqrt{\gamma_{\alpha,\beta}}}+\pmb{\eta}_n}{\frac{1}{\sqrt{\gamma_{\alpha,\beta}}}-{\Big(\frac{1}{n}\log\big(\frac{2}{\delta}\big)\Big)}^\frac{1}{4}\frac{\|\pi_N f\|_\infty}{\|\pi_N f\|_\omega}}
\end{equation}

\end{theorem}

\noindent
{\bf Proof:}
Since for any $x\in I,$ $\pi_N f(x)=\sum_{k=1}^{N+1} c_k(f) \wJ_{k-1}(x),$ where $c_k(f)= \int_I f(t) \wJ_{k-1}(t)\omega_{\alpha,\beta}(t)\, dt,$ then we have 
\begin{equation}\label{system1}
B_N \cdot  \pmb C = \Big[\frac{1}{\sqrt{n}} \wJ_{k-1}(X_j) \Big]_{j,k} \cdot \pmb C =\mathbf P.
\end{equation}
where $\pmb C= \big(c_0(f),\ldots, c_{N}(f)\big)',$ and $\mathbf P = \frac{1}{\sqrt{n}} \big(\pi_N f(X_1),\ldots,\pi_N f(X_n)\big)'.$
The estimator $\widehat f_{n,N}$ is given by ${\displaystyle \widehat f_{n,N} (x)=\sum_{k=1}^{N+1} \widehat c_k \wJ_{k-1}(x),\, x\in I}$ where $\widehat{\pmb C} = \big(\widehat c_0,\ldots, \widehat c_{N}\big)'$ is the least square solution of the over-determined system 
\begin{equation}\label{perturbed_system1}
B_N \cdot \widehat{\pmb C} = \Big[\frac{1}{\sqrt{n}} \wJ_{k-1}(X_j) \Big]_{j,k} \cdot  \widehat{\pmb C} =\mathbf Z_n.
\end{equation}
This can  be written as 
\begin{equation}
B_N \cdot \widehat{\pmb C} =\mathbf P+ \mathbf {\Delta P}.
\end{equation}
where $\mathbf {\Delta P}=\frac{1}{\sqrt{n}} \big((f-\pi_N f)(X_1)+\varepsilon_1,\ldots,(f-\pi_N f)(X_n)+\varepsilon_n\big)'.$
Consequently, the least square solution of \eqref{perturbed_system1} is a perturbation of the least square solution of \eqref{system1}.
In this case, we have (see for example \cite{Zi}) 
$$\frac{\|\pmb C-\widehat{\pmb C}\|_{\ell_2}^2}{\|\pmb C \|_{\ell_2}^2}\leq \kappa_2 (A_N) \frac{\|\mathbf {\Delta P}\|_{\ell_2}^2}{\|\mathbf P\|_{\ell_2}^2}.$$
On the other hand, since $${(f-\pi_N f)(X_j)+\varepsilon_j}^2 \leq 2[{(f-\pi_N f)(X_j)}^2+{\varepsilon_j}^2] ,\quad \forall \; 1\leq j \leq n,$$
then  we have
$$\|\mathbf {\Delta P}\|_{\ell_2}^2\leq 2(\|\mathbf {\widetilde{\Delta} P}\|_{\ell_2}^2+\pmb{\eta}_n^2)$$ 
where $\pmb\eta_n=\displaystyle\max_{\substack{1\leq j\leq n}}|\varepsilon_j|$ and
 $\mathbf{{\widetilde{\Delta} P}}=\frac{1}{\sqrt{n}} \big((f-\pi_N f)(X_1),\ldots,(f-\pi_N f)(X_n)\big)'.$
Hence, one gets 	
\begin{equation}\label{l2norm}	
\frac{\|\pmb C-\widehat{\pmb C}\|_{\ell_2}^2}{\| \pmb C \|_{\ell_2}^2}\leq 2\kappa_2 (A_N) \frac{\|\mathbf {\widetilde{\Delta} P}\|_{\ell_2}^2+\pmb{\eta}_n^2}{\|\mathbf P\|_{\ell_2}^2}.
\end{equation}
In order to estimate $\| f-\widehat f_{n,N} \|_\omega,$ one needs to estimate $\|\pmb C-\widehat{\pmb C}\|_{\ell_2}^2$ which requires the estimation of $\|\mathbf {\widetilde{\Delta} P}\|_{\ell_2}^2+\pmb{\eta}_n^2$ and ${\|\mathbf P\|_{\ell_2}^2}$.
Taking $g=f-\pi_N f$ in Lemma \ref{Lemme1} gives, with probability at least 
$(1-\delta)^2,$
$$\left|\frac{\gamma_{\alpha,\beta}}{n} \sum_{j=1}^n \big((f-\pi_N f)(X_j)\big)^2-\| f-\pi_N f\|_\omega^2\right|\leq \frac{\gamma_{\alpha,\beta}}{\sqrt{n}}  \sqrt{\log\big(\frac{2}{\delta}\big)} \|f-\pi_N f\|_\infty^2.$$
Hence, with probability at least $(1-\delta)^2,$ we have
 $$\frac{1}{n} \sum_{j=1}^n \big((f-\pi_N f)(X_j)\big)^2 \leq\sqrt{\frac{\log\big(\frac{2}{\delta}\big)}{n}}\|f-\pi_N f\|_\infty^2 +\frac{\| f-\pi_N f\|_\omega^2}{\gamma_{\alpha,\beta}}.$$ 
 Equivalently, we have
 $$\|\mathbf {\widetilde{\Delta} P}\|_{\ell_2}^2+\pmb{\eta}_n^2\leq \sqrt{\frac{\log\big(\frac{2}{\delta}\big)}{n}}\|f-\pi_N f\|_\infty^2 +\frac{\| f-\pi_N f\|_\omega^2}{\gamma_{\alpha,\beta}}+\pmb{\eta}_n^2$$ with probability at least 
 $(1-\delta)^2.$  In the same way, taking $g=\pi_N f$ in Lemma \ref{Lemme1}, one gets with probability at least 
 $(1-\delta)^2,$ 
 $$\frac{1}{n} \sum_{j=1}^n \big((\pi_N f)(X_j)\big)^2\geq\frac{\|\pi_N f\|_\omega^2}{\gamma_{\alpha,\beta}}-\sqrt{\frac{\log\big(\frac{2}{\delta}\big)}{n}}\|\pi_N f\|_\infty^2$$ 
 or equivalently 
 $${\|\mathbf P\|_{\ell_2}^2\geq\frac{\|\pi_N f\|_\omega^2}{\gamma_{\alpha,\beta}}-\sqrt{\frac{\log\big(\frac{2}{\delta}\big)}{n}}\|\pi_N f\|_\infty^2}.$$
 Since by Parseval's equality, we have  $\|\pmb C \|_{\ell_2}^2=\|\pi_N f\|_\omega^2 $ and by using \eqref{l2norm}, one gets 
 $$\|\pmb C-\widehat{\pmb C}\|_{\ell_2}^2\leq 2\kappa_2 (A_N) \frac{\sqrt{\frac{\log\big(\frac{2}{\delta}\big)}{n}}\|f-\pi_N f\|_\infty^2 +\frac{\| f-\pi_N f\|_\omega^2}{\gamma_{\alpha,\beta}}+\pmb{\eta}_n^2}{\frac{1}{\gamma_{\alpha,\beta}}-\sqrt{\frac{\log\big(\frac{2}{\delta}\big)}{n}}\frac{\|\pi_N f\|_\infty^2}{\|\pi_N f\|_\omega^2}}$$ with probability at least 
 $(1-\delta)^2.$  To conclude for the proof, it suffices to note that 
 $$\| f-\widehat f_{n,N} \|_\omega \leq \|f-\pi_N f\|_\omega+\|\pi_N f-\widehat f_{n,N}\|_\omega$$
 and use the fact  that $\|\pi_N f-\widehat f_{n,N}\|_\omega=\|\pmb C-\widehat{\pmb C}\|_{\ell_2}.$ \\
 
 The following corollary provides us with more explicit estimation error for the estimator $\widehat f_{n,N}$ in the special
 case where the regression function has a Lipschitzian $p-$th derivative.
 
 \begin{corollary}
 Assume that  $f$ has $p$ continuous derivatives on $I$ and its $p-$th derivative  $f^{(p)} \in Lip(\gamma), 0<\gamma<1.$  Assume that 
 $\|f\|_\omega$ and $\|f\|_\infty$ are bounded away from zero. Moreover, we assume that  $p+\gamma \geq \max(\mu+\frac{1}{2},\frac{1}{2}-\tau),\, \mu=\max(\alpha,\beta),\, \tau=\min(\alpha,\beta).$ Then under the previous notation and assumption,  for any $0<\delta <1$ and for sufficiently large values of 
 $n, N,$ we have with high probability,
 \begin{equation}
 \label{error2}
 \| f-\widehat f_{n,N} \|_\omega \lesssim \frac{\log N}{N^{p+\gamma}}+ \sqrt{2 \kappa_2(A_N)} \frac{\frac{\log N}{N^{p+\gamma-\mu-\frac{1}{2}}} +\pmb{\eta}_n}{\frac{1}{\sqrt{\gamma_{\alpha,\beta}}}-2{\Big(\frac{1}{n}\log\big(\frac{2}{\delta}\big)\Big)}^\frac{1}{4}\frac{\|f\|_\infty}{\|f\|_\omega}}.
 \end{equation} 
  \end{corollary}
 
 \noindent
{\bf Proof:} We recall that under the hypotheses of the corollary, we have see \cite{Prasad},  
	$$\|f-\pi_N f\|_\infty\leq c_1\frac{\log N}{N^{p+\gamma-\mu-\frac{1}{2}}}, \quad \| f-\pi_N f\|_\omega\leq c_2\frac{\log N}{N^{p+\gamma}}.$$
 Here,  $c_1, c_2$ are two positive constants. Hence, for sufficiently large values of $n, N$ and by using the previous inequalities together with
 inequality \eqref{Error1}, one gets the desired error estimate \eqref{error2}.\\
 
 Next, assume that for some $c>0,$ the regression function is the restriction to $I$ of a $c-$bandlimited function which we denote also by $\overline{f}.$ 
 That is $\overline{f}\in L^2(\mathbb R)$ with Fourier transform supported on the compact interval $[-c,c].$ Then a more explicit estimation error of 
 the estimator $\widehat f_{n,N}$ is given in this case by the following corollary.
 
 \begin{corollary} Assume that for some real $c>0,$ $f$ is the restriction to $I$ of a $c-$bandlimited function $\overline{f}.$ Then, for $\beta=\alpha\geq -\frac{1}{2},$ any $0<\delta <1$ and for sufficiently large values of   $n, N,$ we have with high probability,
 \begin{equation}
 \label{error3}
 \| f-\widehat f_{n,N} \|_\omega \lesssim \eta_N^{(\alpha)}\|\overline{f}{f}\|_\omega+ \sqrt{2 \kappa_2(A_N)} \frac{{\Big(\frac{1}{n}\log\big(\frac{2}{\delta}\big)\Big)}^\frac{1}{4} +\frac{\eta_N^{(\alpha)}\|\tilde{f}\|_\omega}{\sqrt{\gamma_{\alpha,\beta}}}+\pmb{\eta}_n}{\frac{1}{\sqrt{\gamma_{\alpha,\beta}}}-2{\Big(\frac{1}{n}\log\big(\frac{2}{\delta}\big)\Big)}^\frac{1}{4}\frac{\|f\|_\infty}{\|f\|_\omega}},
 \end{equation}
where $\eta_N^{(\alpha)}= c^{-1/2} C_\alpha \Big(\frac{ec}{2N+2}\Big)^{N+2},$ for some constant $C_\alpha.$
 \end{corollary}

 \noindent
 {\bf Proof:} 
 In \cite{Jaming_Karoui_Spektor}, it has been shown that for such a function  $f,$ we have 
	$$\| f-\pi_N f\|_\omega\leq \eta_N^{(\alpha)}\|\overline{f} \|_{L^2(\mathbb R)},\quad 
 	\|f-\pi_N f\|_\infty\leq \gamma_N^{(\alpha)}\|\overline{f} \|_{L^2(\mathbb R)}, $$ 
where  
$\eta_N^{(\alpha)}= c^{-1/2} C_\alpha \Big(\frac{ec}{2N+2}\Big)^{N+2},\, \gamma_N^{(\alpha)}= c^{\alpha} C_\alpha \Big(\frac{ec}{2N+2}\Big)^{N+3/2-\alpha},$
for some constant $C_\alpha.$ Hence, by combining the previous two inequalities and \eqref{Error1}, one gets the desired error estimate \eqref{error3}.\\

\begin{remark}
 Although the quantity ${\displaystyle \pmb \eta_n=\max_{1\leq i\leq n} |\eta_i|},$ given in Theorem 4 and the previous two corollaries 
 depends on $n,$ in practice  and with high probability, the different $\eta_i$ are uniformly bounded independently of $n.$   
This is the case for example when the $\eta_i$ are i.i.d. copies of 
 the largely used centered Gaussian random variable  with variance $\sigma^2.$ In this case, 
for  any fixed $k_0>0$ and for any $i\in \mathbb N,$  we have $|\eta_i|\leq k_0 \sigma$ with probability at least ${\displaystyle 1- \mbox{erf}\Big(\frac{k_0}{\sqrt{2}}\Big)\approx 1- \frac{e^{-k_0^2/2}}{k_0 \sqrt{\pi/2}}}.$ This last quantity is very close to $1$ even for reasonable small values of $k_0.$ Consequently, the results of  Theorem 4 and its two corollaries still hold with high probability, with  the quantity $\pmb \eta_n$  replaced by a uniform constant $\pmb \eta.$  
\end{remark}
Next, in order to have an $L_2$-risk error, we need to define a truncated version of our estimator under some conditions.
Let $r>0$ and fix a constant $0<c<0.63$ and consider two positive integers $n,N$ satisfying the inequality 
\begin{equation}\label {1'}
0.63-m_{\alpha,\beta}^2 \frac{(N+1)^{2\mu+2}\, \log(N+1)}{n}-n^{-r}\geq c>0
\end{equation}
Also, assume that the regression function is almost everywhere bounded by a constant $M$, that is 
\begin{equation}\label{2'}
|f(x)|\leq M, \quad \mbox{a.e.} x \in I.
\end{equation} 
Let $\tilde{f}_{N,M}$ be the truncated version of the estimate $\widehat f_{n,N}$ given by 
\begin{equation} \label{3'}
\tilde{f}_{N,M}(x)=\mbox{sign}(\widehat f_{n,N}(x))\min(M, |\widehat f_{n,N}(x)|)
\end{equation}
Under the usual assumption that the added random noise $(\varepsilon_i)_{1\leq i \leq n}$ are i.i.d. real-valued centered  random variables with variance $\sigma^2$ and independent from the $X_i,$ we have the following theorem that provides us with an estimate of the $L_2$-risk error of the estimator $\tilde{f}_{N,M}$. The proof of this theorem is partly inspired from the techniques developed in \cite{Cohen}.
\begin{theorem}\label{thm5}
	Under the previous notation and hypotheses, we have
	\begin{equation}\label{4'}
	\mathbb E\Big[\| f-\tilde{f}_{N,M} \|_\omega^2\Big]\leq \frac{1}{\gamma_{\alpha,\beta} c^2}\left(\frac{N}{n}\sigma^2+\frac{m_{\alpha,\beta}^2 (N+1)^{2\mu+2}}{n}\|f-\pi_Nf\|_\omega^2\right)+\|f-\pi_Nf\|_\omega^2+4M^2 \gamma_{\alpha,\beta} n^{-r}.
	\end{equation}
\end{theorem}
\noindent
{\bf Proof:} Recall that from \eqref{Ineq2.13}, we have for any $\varepsilon>0$
$$\lambda_{\min}(A_N)\geq \big(0.63- m_{\alpha,\beta}^2\frac{  (N+1)^{2\mu+2}}{n} \log(N+1) \big)-\varepsilon$$ 
with probability at least ${\displaystyle 1-\exp\Big(\frac{-n\varepsilon^2}{2m_{\alpha,\beta}^4(N+1)^{4\mu+4}}\Big)}$. In particular for $\varepsilon=\varepsilon_{r,n,N}=\sqrt{\frac{2r\log n}{n}}(N+1)^{\mu+1}$, and under condition \eqref{1'}, one gets 
\begin{equation}\label{5'}
\mathbb P(\lambda_{\min}(A_N)\geq c)\geq 1-n^{-r}
\end{equation}
As it is done in \cite{Cohen}, we consider the two measurable subsets of $I^n$ denoted by $\Omega_+$ and $\Omega_-$, where $\Omega_+$ is the set of all possible draw $(X_1,\cdots, X_n)$ giving $\lambda_{\min}(A_N)\geq c$, while $\Omega_-$ is the set of possible draw giving $\lambda_{\min}(A_N)< c$. Let $d\pmb \rho$ be the probability measure on $I^n$, given by the tensor product of unidimensional probability measure associated to the $B(\alpha+1,\beta+1)$ distribution, that it is 
\begin{equation}
\label{5-5'}
 d \pmb \rho= \otimes^n d \rho(x_i),\qquad d\rho(x_i)= \frac{1}{\gamma_{\alpha,\beta}} \omega_{\alpha,\beta}(x_i) \mathbf 1_I(x_i)\, dx_i,\qquad \gamma_{\alpha,\beta}=B(\alpha+1,\beta+1) 2^{\alpha+\beta+1}.
\end{equation}
Then, we have 
\begin{equation}\label{6'}
\int_{\Omega_{-}}d\pmb \rho=\mathbb P\big\{(X_1,\cdots, X_n)\in I^n; \lambda_{\min}(A_N)< c \big\}\leq n^{-r}.
\end{equation}
Next, by using \eqref{3'}, it is easy to see that the truncated estimator $\tilde{f}_{N,M}$ satisfies 
\begin{equation}\label{7'}
|f(x)-\tilde{f}_{N,M}|\leq |f(x)-\widehat f_{n,N}(x)|, \quad \forall \,x \in I
\end{equation} 
and 
\begin{equation}\label{8'}
|f(x)-\tilde{f}_{N,M}(x)|\leq|f(x)|+|\tilde{f}_{N,M}|\leq 2M, \quad \forall \,x \in I.
\end{equation}
Hence, we have 
\begin{equation}\label{9'}
\mathbb E\big(\| f-\tilde{f}_{N,M} \|_\omega^2\big)=\int_{\Omega_{+}}\| f-\tilde{f}_{N,M} \|_\omega^2d\pmb \rho+\int_{\Omega_{-}}\| f-\tilde{f}_{N,M} \|_\omega^2d\pmb \rho
\end{equation}
By using \eqref{8'}, one gets 
\begin{equation}\label{10'}
\int_{\Omega_{-}}\| f-\tilde{f}_{N,M} \|_\omega^2d\pmb \rho\leq 4M^2 \gamma_{\alpha,\beta} n^{-r}
\end{equation}
To bound the first quantity of the right hand side of \eqref{9'}, we use \eqref{7'} to get 
\begin{eqnarray}\label{8''}
	\int_{\Omega_{+}}\| f-\tilde{f}_{N,M} \|_\omega^2d\pmb \rho&\leq&\int_{\Omega_{+}}\| f-\widehat f_{n,N} \|_\omega^2d\pmb \rho \nonumber\\
	&\leq &\int_{\Omega_{+}}\| f-\pi_Nf \|_\omega^2d\pmb \rho +\int_{\Omega_{+}}\| \pi_Nf-\widehat f_{n,N} \|_\omega^2d\pmb \rho
\end{eqnarray}
By using Parseval's equality and our previous notation, we have, on $\Omega_{+}$
\begin{eqnarray*}
	\| \pi_Nf-\widehat f_{n,N} \|_\omega^2&=&\|\pmb C-\widehat{\pmb C}\|_{\ell_2}^2
	\leq \|A_N^{-1}\|_2^2\|B_N'\mathbf {\Delta P}\|_{\ell_2}^2\\
	&\leq &\frac{1}{c^2}\|B_N'\mathbf {\Delta P}\|_{\ell_2}^2
\end{eqnarray*}
That is $$\mathbb{E}\Big[\| \pi_Nf-\widehat f_{n,N} \|_\omega^2\Big]\leq\frac{1}{c^2}\mathbb E\Big[\|B_N'\mathbf {\Delta P}\|_{\ell_2}^2\Big].$$
But $$\|B_N'\mathbf {\Delta P}\|_{\ell_2}^2=\frac{1}{n^2}\sum_{k=1}^{N+1}\sum_{j,l=1}^{n}\wJ_{k-1}
(X_j)\Big(g_N(X_j)+\varepsilon_j\Big)\wJ_{k-1}(X_l)\Big(g_N(X_l)+\varepsilon_l\Big)$$
where $g_N=(f-\pi_Nf)\quad\bot\quad\wJ_{k-1}, \quad 1\leq k\leq N+1.$
Using the independence of the $X_j$'s as well as the fact that the $\varepsilon_j$'s are independent of the $X_j$'s, together with the facts that $\mathbb E(\varepsilon_j)=0$ and $\mathbb E\big[\varepsilon_j^2\big]=\sigma^2$, one gets 
\begin{equation}\label{11'}
\mathbb E\Big[\|B_N'\mathbf {\Delta P}\|_{\ell_2}^2\Big]=\frac{1}{n^2}\sum_{k=1}^{N+1}\sum_{j=1}^{n}\mathbb E\Big[\varepsilon_j^2\big(\wJ_{k-1}(X_j)\big)^2\Big]+\frac{1}{n^2}\sum_{k=1}^{N+1}\mathbb E\Big[\sum_{j=1}^{n}\big(\wJ_{k-1}(X_j)\big)^2\big(g_N(X_j)\big)^2\Big].
\end{equation} 
Since ${\displaystyle \mathbb E\Big[\varepsilon_j^2\big(\wJ_{k-1}(X_j)\big)^2\Big]=\mathbb E\big[\varepsilon_j^2\big]\mathbb E\Big[\big(\wJ_{k-1}(X_j)\big)^2\Big]=\sigma^2 \frac{1}{\gamma_{\alpha,\beta}}}$ and since ${\displaystyle \mathbb E\Big[(g_N(X_j))^2\Big]=\frac{1}{\gamma_{\alpha,\beta}}\|g_N\|_\omega^2}$, together with the fact that ${\displaystyle \sum_{k=1}^{N+1}(\wJ_{k-1}(X_j))^2\leq m_{\alpha,\beta}^2 (N+1)^{2\mu+2}}$, then we have 
\begin{equation}
\frac{1}{n^2}\sum_{k=1}^{N+1}\mathbb E\Big[\sum_{j=1}^{n}\big(\wJ_{k-1}(X_j)\big)^2\big(g_N(X_j)\big)^2\Big]\leq \frac{m_{\alpha,\beta}^2 (N+1)^{2\mu+2}}{n \gamma_{\alpha,\beta}}\|g_N\|_\omega^2.
\end{equation}
Finally, by combining the previous equalities and inequalities, together with \eqref{9'}--\eqref{11'}, one gets the desired $L_2-$risk error of the estimator $\tilde{f}_{N,M}.$

Next, we briefly describe how our estimator $\widehat f_{n,N}$ can be used in the framework of a general random sampling set $\{X_i,\, 1\leq i\leq n\},$ drawn from a distribution law, other than the Beta law with parameters $\alpha+1,\beta+1.$ There is two cases to consider. In the first case, we assume that the $X_i$ are i.i.d. copies of a random variable with known cumulative distribution function (CDF) $F_X(\cdot).$ In this case, it is well known that the transformed random variables $x_i= F_X(X_i)$ are i.i.d. copies of a random variables following the uniform law over $(0,1).$ Moreover, it is well known that the CDF associated with the Beta distribution with parameters $\alpha+1,\beta+1$ is given by $I_x(\alpha+1,\beta+1):$ the regularized incomplete Beta function. Since this last function is invertible, then the transformed random sampling points, given by 
\begin{equation}\label{newsampling1}
\tau_i = I_x^{-1}(\alpha+1,\beta+1) \big(F_X(X_i)\big),\quad 1\leq i\leq n
\end{equation}
follow the $Beta(\alpha+1,\beta+1)$ distribution. For the interesting special values of $\alpha=\beta=-\frac{1}{2},$ see Remark~\ref{rem2}, the inverse of the $I_x(1/2,1/2)$ is simply given by 
\begin{equation}\label{newsampling2}
I_x^{-1}(1/2,1/2)(t)= \frac{1}{2}\Big(1+\sin\big(\pi t -\frac{\pi}{2}\big)\Big),\quad t\in [0,1].
\end{equation}
In the second case where the CDF $F_X$ is not known, this is the case where the sampling distribution law is unknown,  then one can replace in formula 
\eqref{newsampling1}, the CDF $F_X$ by an estimator of this later. There is a fairly rich literature devoted to build estimators of the CDF of unknown probability laws, see for example \cite{Funke} and the references therein.  
Note that to get accurate  approximate values of the outputs $\widetilde Y_i$ at the transformed random sampling points $\tau_i,$ one may use an interpolation scheme based on the observed $Y_i$ at the original random sampling points $X_i.$ In example 1 of the numerical simulations section, we provide some numerical results that illustrate the stability of the estimator $\widehat f_{n,N}$ when the random sampling set is drawn from the standard normal distribution. \\

Finally, we briefly describe our random pseudoinverse scheme based estimator  can be  generalized  to the multivariate case, with random sampling sets $\{X_i, 1\leq i\leq n\}\subset \mathbb R^d.$ In fact, we only need to replace  
each Jacobi polynomial $\wJ,$ $0\leq k\leq N$  by its tensor product $d-$dimensional version 
\begin{equation*}
\Psi_{\pmb m}^{\alpha,\beta}(\pmb x)= \prod_{j=1}^d \widetilde P^{(\alpha,\beta)}_{m_j}(x_j),\quad \pmb x=(x_1,\ldots,x_d)\in I^d,\quad \pmb m=(m_1,\ldots,m_d)\in \{0,1,\ldots,N\}^d.
\end{equation*}
Note that  $\{\Psi_{\pmb m}^{\alpha,\beta},\, \pmb m\in \mathbb N_0^d\}$ is  an orthonormal basis of $L^2(I^d,\pmb \omega_{\alpha,\beta}),$ where ${\displaystyle \pmb \omega_{\alpha,\beta}(\pmb x) =\prod_{j=1}^d  \omega_{\alpha,\beta}(x_j)}.$ For the sake of simplicity,
we restricted ourselves to the case $d=1.$ The extension of this scheme to higher dimensions is the subject of a future work.

\section{Random pseudo-inverse based scheme for functional regression estimator}

In this paragraph, we first recall that for a compact interval $J,$  the LFR problem  is given by 
\begin{equation}\label{problem2-2}
  Y_i =\int_J X_i(s)\,  \beta_0(s)\, ds +\varepsilon_i,\quad i=1,\ldots,n,
 \end{equation}
where  the $X_i(\cdot)\in L^2(J)$ are random functional predictors, given by ${\displaystyle X_i(s)= \sum_{k=1}^\infty \xi_k Z_{i,k} \varphi_k(s).}$ Here, 
the ${\displaystyle \varphi_k(\cdot)}$ is an orthonormal family of $L^2(J),$
 the   $Z_{i,k}$ are i.i.d centered  random variables with variance $\sigma_Z^2$ and $(\xi_k)_{1\leq k\leq N}$ is a real valued deterministic  sequence.
The $\varepsilon_i$ are i.i.d centered white noise independent  of the $X_i(\cdot)$ and 
$\beta_0(\cdot)\in L^2(J)$ is the unknown slope function to be recovered.  We use a partition of the set $[[1,N]]$ given by $K_N=\Big[\frac{\log N}{\log 2}\Big]+1$ subsets   $ I_k=[[2^{k-1}+1,\min(2^k,N)]].$  On each subset 
	$I_k,$ we consider the random matrix $G_k$ given by 
\begin{equation}
\label{MatrixGk}
G_k=\displaystyle\left[\frac{1}{n} \sum_{i=1}^n \xi_j\xi_l Z_{i,j}Z_{i,l} \right]_{j,l\in I_k}
\end{equation}	
Then, we have 
\begin{equation}
\label{Estimator2}
\widehat \beta_{n,N}(s)=\sum_{k=1}^{K_N} \sum_{j\in I_k} \widehat c_{n,j} \varphi_j(s),\quad \widehat{\pmb c}_{n,k} =[\widehat c_{n,j}]_{j\in I_k}', 
\end{equation}
where
\begin{equation}
\label{Estimator2_coeff}
 \widehat{\mathbf{c}}_{n,k}=G_k^{-1}(F_k'\frac{1}{\sqrt{n}}[Y_i^k]'_{1\leq i\leq n}),\quad 1\leq k \leq K_N.
\end{equation}
Here,  the   $\pmb Y^k = [Y_i^k]_{1\leq i\leq n}$ are obtained from  $\pmb Y=[Y_i]_{1\leq i\leq n}$ by substituting in \eqref{problem2},  $X_i(s)$ with its projection ${X_i^k(s)=\sum_{j\in I_k} \xi_j Z_{i,j} \varphi_j(s)}.$  Consequently, the model \eqref{problem2-2} is substituted with the following model
\begin{equation}\label{Problem2-2}
  Y^k_i =\int_J X^k_i(s)\,  \beta_0(s)\, ds +\varepsilon^k_i,\quad i=1,\ldots,n,\quad 1\leq k\leq K_N,
 \end{equation}
 where the i.i.d. noises $\varepsilon^k_i$  are centered and $\mathbb E\Big[\big(\varepsilon^k_i\big)^2\Big]=\sigma_k^2.$ Moreover, we assume that the $\varepsilon^k_i$ are independent of the $Z_{ij}.$
 Then we show that the estimator $\widehat \beta_{n,N}(\cdot),$ given by 
\eqref{Estimator2} and \eqref{Estimator2_coeff} is an accurate and stable estimator for the slope function $\beta_0(\cdot),$ solution of the previous LFR problem.
The following theorem provides us with a bound for the $2-$condition number of $G_k,$  as well as a bound for the estimation error 
	$\| \beta_0(\cdot)-\widehat \beta_{n,N}(\cdot)\|^2_2.$ For this purpose, we use the function  ${\displaystyle \pmb \xi(\cdot)=\sum_{j=1}^N \xi_j \varphi_j(\cdot)}.$ Moreover, for $1\leq k\leq K_N,$ we use the notation $g^k(\cdot)=\pi_{I_k} g(\cdot),$ the orthogonal projection of $g(\cdot)$ over $\mbox{Span}\{\varphi_j(\cdot),\, j\in I_k\}.$

	\begin{theorem}
		Assume that ${\displaystyle \max_{1\leq j\leq N}{|Z_{ij}|}\leq M}$ with probability at least $(1-\delta_N).$ Then, for any $\eta>0,$ we have with probability at least $1-2\Big(\exp\big(\frac{-n\eta^2}{2M_{\pmb{\xi}}}\big)+\delta_N\Big) $
\begin{equation}
\label{Ineq4.1}
\kappa_2(G_k)\leq \frac{1.72\max_{\substack{j\in I_k}}\sigma_Z^2\xi_j^2+\frac{M_{\pmb{\xi}}}{n}\log(2^{k-1})+\eta}{0.63\min_{\substack{j\in I_k}}\sigma_Z^2\xi_j^2-\frac{M_{\pmb{\xi}}}{n}\log(2^{k-1})-\eta} ; \quad 1\leq k\leq K_N.
\end{equation}
Here $M_{\pmb{\xi}}=M^2\max_{\substack{1\leq j\leq  N}}{|\xi_j|}\|\pmb \xi\|_{\ell_1}.$ Moreover,
  we have 
  \begin{equation}\label{Ineq4.2}
  \|\widehat{\beta}_{n,N}(\cdot)-\beta_0(\cdot)\|_{2}^2\leq
  \frac{1}{n}\Big(\sum_{k=1}^{K_N}\kappa_2(G_k)\Big)  \max_{1\leq k\leq K_N}\frac{\|\beta^k_0(\cdot)\|_{2}^2
  \|(\varepsilon^k_i)_i\|^2_{\ell_2}}{\sigma_Z^2\sum_{j\in I_k} \xi_j^2 c_j^2-\eta}. 
\end{equation}  
with probability at  least $\Big(1-\exp\big(-2n\frac{\eta^2}{ M^4\|\pmb \xi(\cdot)\|_{\ell_2}^4 \|\beta_0(\cdot)\|_{L_2}^4}\big)\Big)^{K_N}.$
  	\end{theorem}

\noindent
{\bf Proof:} Since the i.i.d. random variables $Z_{i,j}$ are centred with variances $\sigma_Z^2$, then $$\mathbb{E}(G_k)=\begin{bmatrix}\sigma_Z^2\xi_{2^{k-1}+1}^2&& \\&\ddots\\&&\sigma_Z^2\xi_{2^k} \end{bmatrix}.$$ 
Consequently, the minimum and  the maximum eigenvalue of $\mathbb{E}(G_k)$ are given by $$\Lambda_{\min}=\lambda_{\min}(\mathbb{E}(G_k))=\sigma_Z^2\displaystyle\min_{\substack{j\in I_k}}{\xi_j^2}, \quad \Lambda_{\max}=\lambda_{\max}(\mathbb{E}(G_k))=\sigma_Z^2\displaystyle\max_{\substack{j\in I_k}}{\xi_j^2}.$$
On the other hand, the random matrix $G_k$ is written in the following form $$G_k=\frac{1}{n}\sum_{i=1}^{n}D_{i,k}, \quad D_{i,k}=\left[\sum_{i=1}^n \xi_j\xi_l Z_{i,j}Z_{i,l} \right]_{j,l\in I_k}.$$
Note that each matrix $D_{i,k}$ is  positive semi definite. This follows from the fact that for any $\pmb x\in \R^{2^{k-1}}$, we have $$\pmb{x'}D_{i,k}\pmb{x}=\pmb{x'} B_{i,k} B_{i,k}\pmb{x}\geq 0, \quad B_{i,k}=\frac{1}{\sqrt{n}}[\xi_jZ_{i,j}]_{j\in I_k}.$$
By using Gershgorin circle theorem, one gets
$$\begin{array}{lll} 
\lambda_{\max}(D_{i,k})&\leq &\frac{1}{n}\max_{\substack{j\in I_k}}{|\xi_jZ_{i,j}|}\displaystyle\sum_{l=2^{k-1}+1}^{2^k}|\xi_l Z_{i,l}|\\
&\leq&\frac{M^2}{n}\max_{\substack{j\in I_k}}{|\xi_j|}\displaystyle\sum_{l=2^{k-1}+1}^{2^k}|\xi_l|\leq \frac{M_{\pmb{\xi}}}{n}
\end{array}$$
Hence, we have
$$\mathbb{E}(\lambda_{\min}(G_k))\geq 0.63\min_{\substack{j\in I_k}}\sigma_Z^2\xi_j^2-\frac{M_{\pmb{\xi}}}{n}\log(2^{k-1})$$
and $$\mathbb{E}(\lambda_{\max}(G_k))\leq 1.72\max_{\substack{j\in I_k}}\sigma_Z^2\xi_j^2+\frac{M_{\pmb{\xi}}}{n}\log(2^{k-1})$$
Then, we use McDiarmid's concentration inequality to conclude that  for any $\eta>0$, we have 
$$\lambda_{\min}(G_k)\geq \mathbb{E}(\lambda_{\min}(G_k))-\eta,\quad \lambda_{\max}(G_k)\leq \mathbb{E}(\lambda_{\max}(G_k))+\eta$$
with probability at least $1-\exp\frac{-n\eta^2}{2M_{\pmb{\xi}}}$. By combining the previous two inequalities, one gets \eqref{Ineq4.1}.\\

Next to get the estimation error given by \eqref{Ineq4.2}, we note that on each subset $ I_k$, the associated  exact expansion coefficients vector of 
$\beta_0(\cdot)$ is given by  $\mathbf{c}^k=[c_j]_{j\in I_k}.$ Similarly, we define  the expansion coefficients restricted to $I_k$ of the estimator 
$\widehat \beta_{n,N},$ that we denote by  $\widehat {\mathbf{c}}_{n,k}.$ These two coefficients vectors are given by
\begin{equation}\label{Ineq4.3}
{\mathbf{c}^k}'=G_k^{-1}\Big(F_k'\frac{1}{\sqrt{n}}[Y_i^k-\varepsilon_i^k]'_{1\leq i\leq n}\Big)=G_k^{-1}\Big(F_k'\frac{1}{\sqrt{n}}[\tilde{Y_i^k}]'_{1\leq i\leq n}\Big), \quad \widehat{\mathbf{c}}_{n,k}=G_k^{-1}\Big(F_k'\frac{1}{\sqrt{n}}[Y_i^k]'_{1\leq i\leq n}\Big),
\end{equation}
where,
\begin{equation}
\label{Ineq4.4}
F_k=\frac{1}{\sqrt n}\big[\xi_jZ_{i,j}\big]_{1\leq i\leq n ,j\in I_k},\quad 1\leq i\leq n,\quad 1 \leq k\leq K_N.
\end{equation}
 
Then we have the following perturbation result for the pseudo-inverse least square solution of  an over-determined system,
\begin{equation}
\label{Ineq4.5}
\frac{\|\widehat{\mathbf{c}}_{n,k}-\mathbf{c}^k\|_{\ell_2}^2}{\|\mathbf{c}^k\|_{\ell_2}^2}\leq \kappa_2(G_k)\frac{\|\mathbf{\widetilde{Y}}^k-\mathbf{Y}^k\|_{\ell_2}^2}{\|\mathbf{\widetilde{Y}}^k\|_{\ell_2}^2}=\kappa_2(G_k)\frac{\frac{1}{n}\|(\varepsilon_i^k )_i\|_{\ell_2}^2}{\frac{1}{n}\|(\widetilde{Y_i}^k)_i\|_{\ell_2}^2}.
\end{equation}
On the other hand, since $\Big(\widetilde{Y_i}^k\Big)^2=\Big(\sum_{j\in I_k}\xi_j Z_{i,j}c_j\Big)^2,$ then we have 
$$0\leq \big(\widetilde{Y_i}^k\big)^2\leq M^2 \|\pmb \xi^k(\cdot)\|^2_{\ell_2} \|\beta_0^k(\cdot)|_{L_2}^2 $$
and $$\mathbb{E}\big[\big(\widetilde{Y_i}^k\big)^2\big]=\sigma_Z^2 \sum_{j\in I_k} \xi_j^2 c_j^2$$
Hence by using Hoeffding's inequality, one concludes that for any $\eta>0$, we have
\begin{equation}
\label{Ineq4.6}
\frac{1}{n}\sum_{i=1}^{n}\big(\widetilde{Y_i}^k\big)^2-\mathbb{E}\big[\big(\widetilde{Y_i}^k\big)^2\big]=\frac{1}{n}\sum_{i=1}^{n}\big(\widetilde{Y_i}^k\big)^2-\sigma_Z^2 \sum_{j\in I_k} \xi_j^2 c_j^2 \geq -\eta
\end{equation}
with probability at least $1-\exp\big(-2n\frac{\eta^2}{{M}^4\|\xi^l(\cdot)\|^2_{\ell_2} \|\beta_0^l\|_{L_2}^4}\big)$. Finally, by using the inequalities \eqref{Ineq4.3}--\eqref{Ineq4.6}, together with Parseval's equality and some straightforward computations,  one gets the desired inequality \eqref{Ineq4.2}.

\begin{remark}
In practice, the deterministic sequence $\pmb \xi =(\xi_j)_j$ satisfies a decay condition of the type $\xi_j \asymp j^{-s}$ for some $s>0.$ 
In this case,  the combined   dyadic decomposition and random pseudo-inverse scheme estimator $\widehat{\beta}_{n,N}(\cdot)$ is stable 
 in the sense that 
 \begin{equation}
 \label{Ineq4.7}
 \kappa_2(G)= \sum_{k=1}^{K_N}\kappa_2(G_k) \lesssim 2^s \frac{1.72 \cdot \log N}{0.63\cdot  \log 2}  
 \end{equation}
 is a relatively  small cumulative condition number.
\end{remark}
Next, we are interested in having an $L_2$-risk error of a truncated version of our LFR estimator $\widehat \beta_{n,N}(\cdot).$ More precisely, we assume that the true slope function $\beta_0(\cdot)$ is almost every where bounded, that is there exists $L>0$ such that
$$|\beta_0(x)| \leq L,\quad a.e. x\in J.$$
Then, as for the nonparametric estimator of the previous section, the truncated  $\widehat \beta_{n,N}(\cdot)$ estimator, denoted by $\widetilde \beta_{N,L}(\cdot) $ is defined by 
$$ \widetilde \beta_{N,L}(x)=\mbox{Sign}\big(\widehat \beta_{n,N}(x)\big )
\min\Big(L, |\widehat \beta_{n,N}(x)|\Big).$$
Let $r > 0$ and $\eta_k >0$ be such that 
\begin{equation}\label{Ineqq4.7}
\mathbb{P}\big(\lambda_{\min}(G_k)\geq \eta_k\big)\geq 1-n^{-r},\quad 1\leq k\leq K_N.
\end{equation}
The following theorem provides us with the  $L_2$-risk of the estimator $\widetilde{\beta}_{N,M}.$ 
\begin{theorem}
	Under the previous notations and hypotheses, we have
	\begin{equation}
	\mathbb E\Big[\| \widetilde{\beta}_{N,M}(\cdot)-\beta_0(\cdot)\|_{2}^2\Big]\leq \frac{\sigma_Z^2\|(\xi_i)_i\|_{\ell_2}^2}{n^2}\Bigg(\sum_{k=1}^{K_N}\frac{\sigma_k^2}{\eta_k^2}|I_k|\Bigg)+\frac{4 L^2}{n^r}K_N.
	\end{equation}
\end{theorem}
\noindent
{\bf Proof:} Let $\mathbf{\tilde{c}_{n,k}}$ be the expansion coefficients vector of $ \widetilde{\beta}_{N,M}^k$, the projection over $\mbox{Span}\{\varphi_j(\cdot),\, j\in I_k\}$ of $ \widetilde{\beta}_{N,M}.$  For each $1\leq k\leq K_N,$ let $\Omega_{+,k}$ and $\Omega_{-,k}$  be the set of all possible draw $(X_1(\cdot),\cdots, X_n(\cdot))$ for which $\lambda_{\min}(G_k)\geq \eta_k$  and $0<\lambda_{\min}(G_k)< \eta_k$ respectively. Then , we have
\begin{equation}\label{1"}
\mathbb E\Big[\|\mathbf{\tilde{c}_{n,k}}-\mathbf{c}^k\|_{\ell_2}^2\Big]\leq \int_{\Omega_{+,k}}\int_J|( \widetilde{\beta}_{N,M}^k-\beta_0^k)(x)|^2dx d\pmb \rho +\int_{\Omega_{-,k}}\int_J|( \widetilde{\beta}_{N,M}^k-\beta_0^k)(x)|^2dx d\pmb \rho.
\end{equation}
where 
\begin{equation}\label{2"}
\int_{\Omega_{-,k}}\int_J|( \widetilde{\beta}_{N,M}^k-\beta_0^k)(x)|^2dx d\pmb \rho \leq \frac{4 L^2}{n^r}.
\end{equation}
and 
\begin{equation}\label{3"}
\int_{\Omega_{+,k}}\int_J|( \widetilde{\beta}_{N,M}^k-\beta_0^k)(x)|^2dx d\pmb \rho \leq \int_{\Omega_{+,k}}\int_J|(\widehat{\beta}_{n,N}^k-\beta_0^k)(x)|^2dx d\pmb \rho.
\end{equation}
Here, $d\pmb \rho$ is the tensor product  probability measure, defined in a similar manner as \eqref{5-5'} and  associated 
with the probability law of the random predictors $X_i(\cdot).$ 
By using Parseval's equality, we have from \eqref{Ineq4.3}
$$\int_J|(\widehat{\beta}_{n,N}^k-\beta_0^k)(x)|^2dx=\|\widehat{\mathbf{c}}_{n,k}-\mathbf{c}^k\|_{\ell_2}^2\leq \frac{1}{n}\|G_k^{-1}\|_2^2\|F_k'\big(\mathbf{\widetilde{Y}}^k-\mathbf{Y}^k\big)\|_{\ell_2}^2 .$$
So that on $\Omega_{+,k},$ we obtain 
$$\int_J|(\widehat{\beta}_{n,N}^k-\beta_0^k)(x)|^2dx\leq \frac{1}{n \eta_k^2}\|F_k'\big(\mathbf{\widetilde{Y}}^k-\mathbf{Y}^k\big)\|_{\ell_2}^2.$$
But $$F_k'\big(\mathbf{\widetilde{Y}}^k-\mathbf{Y}^k\big)=\frac{1}{\sqrt{n}}\Bigg[\displaystyle\sum_{l=1}^{n}-\varepsilon_l^k\xi_iZ_{l,i}\Bigg]'_{i\in I_k}=\frac{1}{\sqrt{n}}\Big[a_i\Big]'_{i\in I_k}.$$ 
By using the hypotheses on the noises $\varepsilon_i,$ it is easy to see that 
$$\mathbb E\big[a_i^2\big]=\mathbb E\Bigg[ \sum_{l,j=1}^{n}\varepsilon^k_l\varepsilon^k_j\xi_l\xi_iZ_{l,i}Z_{j,i}\Bigg]=\sigma_k^2\sigma_Z^2\sum_{j=1}^{n}\xi_j^2.$$
Consequently, we have
\begin{equation*}
\int_J|(\widehat{\beta}_{n,N}^k-\beta_0^k)(x)|^2dx\leq \frac{|I_k|}{n^2\eta_k^2}\sigma_k^2\sigma_Z^2\|(\xi_i)_i\|_{\ell_2}^2.
\end{equation*}
The previous inequality, together with \eqref{1"}, \eqref{2"} and \eqref{3"} lead to
$$\mathbb E\Big[\|\mathbf{\tilde{c}_{n,k}}-\mathbf{c}^k\|_{\ell_2}^2\Big]\leq \frac{|I_k|}{n^2\eta_k^2}\sigma_k^2\sigma_Z^2\|(\xi_i)_i\|_{\ell_2}^2+\frac{4 L^2}{n^r}$$

Finally,  by using Parseval's equality and by adding inequalities on subsets $I_k$, one gets the desired $L_2-$risk error of the estimator $ \widetilde{\beta}_{N,M}$.
\begin{remark}
	For the special of  $K_N=1,$ that is  our estimator $\widehat \beta_{n,N}(\cdot)$ is constructed without the dyadic partition of the set $[[1,N]],$ then the $L_2-$risk of the previous theorem is simply given by
	$$	\mathbb E\Big[\| \widetilde{\beta}_{N,M}(\cdot)-\beta_0(\cdot)\|_{2}^2\Big]\leq \frac{\sigma^2\sigma_Z^2\|\xi\|_{\ell_2}^2N}{n^2 \eta^2_N}+\frac{4 L^2}{n^r}. $$
\end{remark}

\section{Numerical Simulations}

In this section, we first give some numerical  simulations that illustrate the different results of this paper. Then, we provide real data application of our nonparametric regression estimator with data corresponding to the Covid-19 daily  spread over $15$ months ( starting from March 1st 2020) of four  countries covering different regions of the world. 

\subsection{ Numerical examples}

\noindent
{\bf  Example 1: } In this first example, we illustrate the result of Theorem~3, as well as an adaptation of this result when the random sampling set is drawn from a standard normal distribution, as described by \eqref{newsampling1} and \eqref{newsampling2}.
Then, we illustrate by some numerical test, the inequality \eqref{Ineq4.7} which is a consequence of Theorem~5. For this purpose,  we have  considered the case of the re-normalized  Jacobi polynomials, so that they 
are orthonormal over $I_1=[0,1]$ with   the two special cases of $\alpha=\beta=-\frac{1}{2}$ as well as $\alpha=\beta=0.$
In this case, the corresponding coefficients $m_{\alpha,\beta},$ given by Proposition 1 and remark~1, are 
$m_{-\frac{1}{2},-\frac{1}{2}}=\frac{2}{\pi},\quad m_{0,0}= \sqrt{\frac{2+1/\sqrt{2}}{3}}\approx 0.95.$
Then, we have considered the values of $N=5, 10, 15, 20$ and different values of $n$ the size of the sampling set $\{ X_j, j=1,\ldots,n\}$ following the
beta $B(\alpha,\alpha)-$law. Also, we have computed the average of the condition number $\kappa_2(A_N)$ over $50$ such realizations of $A_N.$
Then, we have repeated  the previous simulations, with the initial random sampling set $\{X_i,\, 1\leq i\leq n\}$ drawn from the standard normal distribution. We have considered the transformed 
random sampling ${\tau_i,\, 1\leq i\leq n},$ associated with  the $Beta(\alpha+1,\beta+1)$ probability law. We denote by $\widetilde A_N,$ the counterparts of the random  matrix $A_N$ when the transformed sampling points $\tau_i$ are used instead of the $X_i.$  The obtained numerical results are given by  Table~1.\\

Next, to illustrate the stability of our proposed LFR estimator $\widehat \beta_{n,N}(\cdot),$ we have computed the cumulative $2-$norm condition numbers of the positive definite random matrices $G_k,$ given by \eqref{MatrixGk}, with different values 
of the integers $N, n$ and the decay rate exponent $s$ of the deterministic sequence $(\xi_j)_j.$ Also, we assume that the i.i.d. random variables $Z_j$  follow the uniform law $U(-\sqrt{3},\sqrt{3}).$ In Table 2, we have listed the corresponding values of the cumulative condition numbers ${\displaystyle \kappa_2(G)=\sum_{k=1}^{K_N}\kappa_2(G_k).}$ It is interesting to note that from these numerical values, the actual values of the $\kappa_2(G)$ are slightly smaller than the theoretical bounds given by the right-hand side of inequality \eqref{Ineq4.7}.

   \begin{center}
 \begin{table}[h]
 \vskip 0.2cm\hspace*{2cm}
 \begin{tabular}{cccccccccc} \hline
 &&&&&&&&&\\
    $\alpha=\beta$ &$N$&$n$&$\kappa_2(A_N)$& $\kappa_2(\widetilde A_N)$&$\alpha=\beta$ &$N$&$n$&$\kappa_2(A_N)$&$\kappa_2(\widetilde A_N)$ \\ 
     &&&&&&&&&\\  \hline
   $-0.5$  &$5$  &$25 $ &$7.74 $&$6.48$&  $0.0$  &$5$  &$ 40 $ &$8.36$ &$7.27$      \\
           &$10$ &$40 $ &$17.57 $&$11.03$&         &$10$ &$100 $ &$19.94$ &$13.36$     \\
           &$15$ &$60 $ &$24.72 $&$22.03$&         &$15$ &$125 $ &$33.84$ &$22.52$   \\
           &$20$ &$100$ &$27.09 $&$12.03$&         &$20$ &$250$ &$133.07$ & $29.82$ \\   \hline
  \end{tabular}
  \caption{Illustration of Theorem 3.}
  \end{table}
 \end{center} 
   \begin{center}
  \begin{table}[h]
  \vskip 0.2cm\hspace*{1cm}
  \begin{tabular}{ccccccccccccc} \hline
     $s$    &$N$    &$n$    &$\kappa_2(G)$& $N$    &$n$     &$\kappa_2(G)$& $N$  &$n$  &$\kappa_2(G)$& $N$ &$n$  &$\kappa_2(G)$   \\   \hline
    $0.75$  &$20$   &$100 $ &$12.05 $      &  $30$  &$ 100 $ &$15.49 $      & $40$ &$100$&$17.63$&$50$&$100$&$19.66$      \\
        $$  &$ $    &$150 $ &$11.39$       &        &$ 150 $ &$12.85 $      & $ $  &$150$&$15.66$& $ $ &$150$&$18.08$      \\
        $$  &$ $    &$200$  &$11.22 $      &        &$200 $  &$12.26 $      & $ $  &$200$&$15.16$& $ $ &$200$&$15.90$     \\
           $1.5  $  &$20$   &$100 $ &$25.59 $      &  $30$  &$ 100 $ &$29.55 $      & $40$ &$100$&$36.23$&$50$&$100$&$41.01$      \\
                $$  &$ $    &$150 $ &$23.78$       &        &$ 150 $ &$28.60 $      & $ $  &$150$&$33.32$& $ $ &$150$&$37.11$      \\
                $$  &$ $    &$200$  &$21.90$      &         &$200 $  &$27.97 $      & $ $  &$200$&$31.39$& $ $ &$200$&$31.39$     \\   \hline
         
   \end{tabular}
   \caption{Illustration of the cumulative condition number given by  \eqref{Ineq4.7}.}
   \end{table}
  \end{center} 

\noindent
{\bf Example 2:} {\it  A synthetic test function.}  In this example, for a given real number $s>0,$  we consider   the 
Weierstrass test function, given by 
\begin{equation}\label{test_function}
W_s(x)= \sum_{k\geq 0} \frac{\cos(2^{k}\pi x)}{2^{ks}},\quad -1\leq
x\leq 1.
\end{equation}
Note that $W_s \in C^{s-\epsilon}(I),\,\forall \epsilon <s.$  Here, $C^r(I)$ is the usual H\"older  space defined over the interval $I$ and having a H\"older smoothness exponent $r>0.$ We have considered the two values of $s=1, 2.$ Note that for $s=1,$ $W_s(\cdot)$ is continuous but  nowhere differentiable, which makes this case hard to handle by classical nonparametric regression estimators.  
Then, we considered the values of $\alpha=\beta=-\frac{1}{2}$ and constructed our  estimator $\widehat f_{n,N}(\cdot),$  based on the first $N+1$ Jacobi polynomials $\widetilde P_k^{(\alpha,\alpha)}$ with $N=10, 20, 30$ and 
 a random sampling set of size $n=100$ and  following the beta $B(\alpha+1,\alpha+1)$ distribution on $I.$ The added centered white noise in \eqref{problem1} 
 has variance $\sigma^2,$ with $\sigma=0.1, \, 0.05.$ Also, We have performed $10$ realizations of the numerical  approximations of the Weierstrass function by our estimator.\\
 
 Moreover, for comparison purpose, we have repeated the above numerical simulations 
 by applying the Kernel Ridge Regression (KRR) estimator, associated with the Sinc-kernel, given by 
 ${\displaystyle K_c(x,y)=\frac{\sin(c(x-y))}{\pi(x-y)},\, x,y\in I,}$ where $c>0$ is a tuning parameter 
 called bandwidth. Note that the popular KRR estimator is briefly described as follows.
 Given  a reproducing kernel Hilbert space (RKHS), associated with a positive-definite kernel $K(\cdot,\cdot)$ and  given a convenient regularization parameter $\lambda >0$, the KRR  estimator of $f$ is given by the solution of the following minimization problem,
 \begin{equation}\label{RKK1}
 \widehat f_\lambda=\arg\min_{f\in \mathcal H} \left\{\frac{1}{n} \sum_{i=1}^n \Big(f(X_i)-Y_i \Big)^2 +\lambda \|f\|_{\mathcal H}^2\right \},
 \end{equation}
 where $\|\cdot\|_{\mathcal H}$ is the usual norm of $\mathcal H.$ Thanks to the representer theorem, the solution of the previous minimization 
 problem is given by
 \begin{equation}\label{RKK2}
  \widehat f_\lambda(x)=\sum_{k=1}^n c_k K(x,X_k),\quad \pmb c=[c_1,\ldots,c_n]'= G_\lambda^{-1} \frac{1}{n} [Y_1,\ldots,Y_n]',
   \end{equation}
 where $G_\lambda = \Big[ \frac{1}{n} K(X_i,X_j)\Big]_{1\leq i,j\leq n} +\lambda I_n$ is the regularized Gram-matrix, with $I_n$ is the $n\times n $ identity matrix. The convenient  value of the regularization parameter is chosen by the well-known cross validation technique.
 We denote by $\widehat f_c,$ the KRR estimator based on the Sinc-kernel with parameter $c>0.$
  The average mean squared estimation errors, associated with the estimators $\widehat f_{n,N}$ and $\widehat f_c$  for this example and for different values of $N$ and $c$ are given in Table 3.  
 From these numerical results, one concludes that both estimators have similar accuracy. Nonetheless, our estimator is much faster than the KRR based estimator since it involves the inversion of a relatively small size random matrix. This is not the case for a KRR based estimator, specially when handling large sampling dataset.  As an example, for a moderate size of the sampling set $n=200,$ and with the values of $N=c=10,$ 
 the computing time required for the  construction of the KRR estimator $\widehat f_c$ is approximately 10 times the computing time for our proposed estimator 
$\widehat f_{n,N}.$ 

\begin{center}
\begin{table}[h]
\vskip 0.2cm\hspace*{1cm}
\begin{tabular}{ccccccccccc} \hline
 $\sigma$& $s$     &$N$&$MSE(\widehat f_{n,N})$  &$c$&  $MSE(\widehat f_{c})$&$s$& $N$& $MSE(\widehat f_{n,N})$ &$c$&  $MSE(\widehat f_{c})$ \\   \hline
$0.1$    & $1.0$   &$10$& $1.53 e-1$             &$10$& $1.60e-1 $ &      $2.0$& $10$  & $1.71 e-3$ &$10 $ & $1.44e-3 $    \\
         &         &$20$& $5.88 e-3$             &$20$& $7.53e-3 $ &      $ $  &$20$   & $2.00 e-3$ &$20 $ & $1.58e-3 $    \\
         &         &$30$& $3.94 e-3$             &$30$& $4.05e-3 $ &      $ $  &$30$   & $2.93 e-3$ &$30 $ & $2.12e-3 $    \\
$0.05$    & $1.0$   &$10$& $1.42 e-1$            &$10$& $1.29e-1 $ &      $2.0$& $10$  & $8.23 e-4$ &$10 $ & $9.31e-4 $    \\
         &          &$20$& $4.86 e-3$            &$20$& $6.31e-3 $ &      $ $  &$20$   & $6.35 e-4$ &$20 $ & $6.98e-4 $    \\
         &          &$30$& $1.82 e-3$            &$30$& $1.85e-3 $ &      $ $  & $30$  & $6.55 e-4$ &$30 $ & $9.34e-4 $    \\  \hline
\end{tabular}
\caption{Mean squared estimation errors for example  2.}
\label{tableau4}
\end{table}
\end{center}

\vskip 0.5cm
\noindent
{\bf Example 3:} In this last example, we illustrate the performance of our LFR estimator 
$\widehat \beta_{n,N}(\cdot).$ We consider the  test LFR problem  that has been proposed in \cite{Hall} and used later on for comparisons purposes in \cite{Shin}. For this test problem,  $J=[0,1]$ and the slope function is given by 
$$\beta_0(s) = \sum_{j=1}^{50} c_j \varphi_j(s) =\sum_{j=1}^{50} 4 \frac{(-1)^{j+1}}{j^2} \varphi_j(s), \quad \varphi_j(s)=\left\{\begin{array}{ll} 1& \mbox{ if } j=1\\ \sqrt{2} \cos(\pi j s)&\mbox{ if } j\geq 1.\end{array}\right.$$
The random predictor functional is given by 
${\displaystyle X(\cdot)=\sum_{j=1}^{50} \xi_j Z_j \varphi_j(\cdot),}$ where the deterministic sequence $(\xi_j)_j$
is given by ${\xi_j=\frac{(-1)^{j+1}}{j^{s/2}},\,\, s\geq 0}$ 
and   the $Z_k$ are i.i.d. random sample following the  uniform law $U(-\sqrt{3},\sqrt{3}).$ The added centred white noise $\varepsilon_i$ 
has variance $\sigma^2$ with $\sigma=0.5.$ Then, we have computed our estimator $\widehat \beta_{n,N}(\cdot)$ according our scheme given by 
\eqref{Estimator2}-\eqref{Estimator2_coeff}, with  values of $N=50,$ $n=100, 200, 300$ and $s=1.5, 2.0, 4.0.$   We have computed the average of the squared prediction and estimation error over $10$ realizations. These errors are respectively given by 
$$E_0=\|\widehat{\beta}_0(\cdot)-\widehat{\beta}_{n,N}(\cdot)\|^2_0=\sum_j \frac{1}{j^s} \big(c_j-\widehat c_j\big)^2,\quad 
E_2=\|\widehat{\beta}_0-\widehat{\beta}_{n,N}\|^2_2=\sum_j \big(c_j-\widehat c_j\big)^2.$$

{\tiny
\begin{center}
\begin{table}[h]
\vskip 0.2cm\hspace*{0.1cm}
\begin{tabular}{cccccccccccc} \hline
 $s$   &$n$&$E_0$&$E_2$        &   $s$  &$n$&$E_0$&$E_2$&   $s$  &$n$&$E_0$&$E_2$\\   \hline
$1.5$  & $100$&$2.62 e-3$& $7.72 e-2$ &$2.0$& $100$    & $1.52 e-3$ & $6.94 e-3$& $4.0$ &$100$  &$2.02 e-3$   & $3.32 e-1$    \\
       & $200$&$1.28 e-3$&$ 3.12 e-2$ &        & $200$   &$ 1.11 e-3$&  $5.17 e-3$&     &$200$  &$3.92 e-4$   & $9.11 e-2$\\
       & $300$&$3.15 e-4$& $1.62 e-2$ &        & $300$   & $3.84 e-4$&  $1.54 e-3$&     &$300$  &$2.51 e-4$   & $1.28 e-2$   \\ \hline
\end{tabular}
\caption{Prediction and estimation errors of the estimator $\widehat{\beta}_{n,N}(\cdot)$ for $N=50$ and different values
of $s$, $n.$}
\label{tableau4}
\end{table}
\end{center}
}
The obtained numerical values of these errors are given by Table 4. These numerical results are coherent with the theoretical results given by Theorem 5. 
Moreover, by comparing our results with those obtained in \cite{Shin} for this same test problem, we conclude 
that our estimator $\widehat \beta_{n,N}(\cdot)$ clearly  outperforms in terms of accuracy and convergence speed,  the estimators proposed in the previous references.

\subsection{A Real data Application: Covid-19 world's countries spread. } 

In this last section, we use the publicly available daily updated dataset of world's countries spread of Covid-19, provided by Our World in Data. This dataset contains the updated data starting from March 1st 2020 till October 16th and corresponding to the daily confirmed cases as well as the daily deaths for different world's countries. The  URL of the dataset is   https://github.com/owid/covid-19-data/tree/master/public/data

\begin{figure}[h]\hspace*{0.5cm}
{\includegraphics[width=16cm,height=9.5cm]{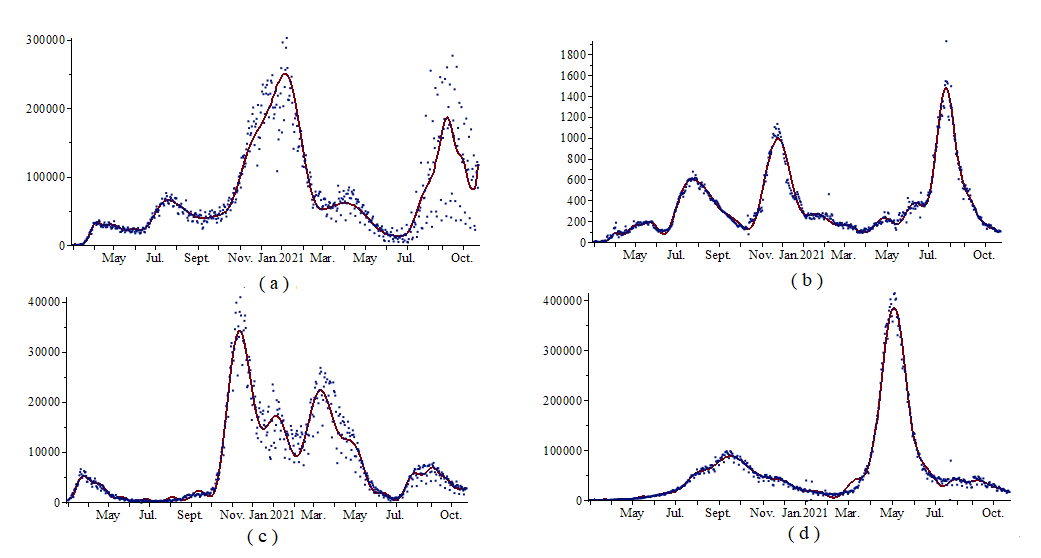}}
 \vskip -0.5cm\hspace*{1cm} \caption{Plots of the true Covid-19 daily cases (points) versus their Jacobi based estimates (line) for (a) USA, (b) Algeria, (c) Italy, (d) India.}
\end{figure} 

For the present simulations, we restrict ourselves to  the data corresponding to daily cases during the period March 1st-October 16th and corresponding to four  countries: USA, Algeria, Italy and  India.  The considered data covers $m=594$ days and the aim is to build a non-parametric regression estimators for these countries daily  covid-19 cases.  To apply our estimator based on  random pseudo-inverse associated with  Jacobi polynomials, we have considered the special values of the parameters $\beta=\alpha=-\frac{1}{2},$ $N=40.$ Then, we have considered a  sampling set ${\displaystyle S=\{ X_i,\, i=1,\ldots,n\}}$ with $n=340$ random samples in the interval $[0,1]$ and following the $B(\alpha+1,\alpha+1)$ law. We found that the average $2-$condition number of the small size $(N+1)\times (N+1)$ random  matrix $A_N = B_N' B_N$ is given by $\kappa_2(A_N)\approx 9.12.$  This ensures the stability 
of the proposed nonparametric regression estimator, as predicted by Theorem 3. We have considered  $n$ observed  outputs $Y_i,$ given by 
$Y_i= f([m X_i]),$ where $f(k)$ is  the observed number of daily  Covid-19 confirmed cases at day $k$ with $1\leq k \leq 594.$ The coefficients $(\widehat c_k)_{0\leq k\leq N}$ of the Jacobi polynomials based  estimators $\widehat f_{n,N}$ are computed in a fast way by
using formula \eqref{coefficients}. Note that since the random sampling points belong to the interval $J=[0,1],$ then our estimators are 
constructed by the use of the re-normalized Jacobi polynomials given by ${\displaystyle Q_k(x)=\frac{1}{\sqrt{2}} \widetilde P_k(2 x-1),\, x\in J.}$ These last set of Jacobi polynomials are orthonormal on $J.$ In order to make this estimator robust (that is not sensitive to data outliers),
we  have applied the RANSAC algorithm trick with $10$ iterations and we have chosen the estimator that best fits the whole data set of 
$594$ measurements from the $n=340$ the size of  random sampling set used at each iterations. In Figure 1, we  have plotted in points, the graphs of the real data corresponding to the daily confirmed cases, versus the graph of the associated Jacobi polynomials based  estimator $\widehat f_{n,N}(\cdot).$   These numerical simulation results are coherent with the theoretical results 
of Theorem 3. In particular, the estimator $\widehat f_{n,N}(\cdot)$ has been made robust by just using few iterations.

\end{document}